\documentclass[12pt]{article}
\usepackage[dvips]{graphicx}
\usepackage{amsmath,amssymb,amsthm}
\usepackage{bm}
\usepackage{latexsym}
\usepackage{cases}

\newtheorem{dfn}{Def}[section]
\newtheorem{thm}[dfn]{Theorem}
\newtheorem{prop}[dfn]{Proposition}
\newtheorem{lem}[dfn]{Lemma}
\newtheorem{rem}[dfn]{Remark}
\newcommand{\for}{\textrm{for}}

\makeatletter

\newcounter{myc}[section]

\newcounter{mycc}[section]

\def\id#1{\def\@id{#1}}
\def\department#1{\def\@department{#1}}
\def\superadvisor#1{\def\@superadvisor{#1}}

\makeatletter
    
    \@addtoreset{equation}{section}
\makeatother

\begin{document}
\title{\large
\bf{Diffusion wave phenomena and $L^p$ decay estimates \\of solutions of compressible viscoelastic system}}
\author{
\large{Yusuke Ishigaki}}
\date{\small
Department of Mathematics, \\
Tokyo Institute of Technology, \\
Tokyo 152-8551, Japan, \\
e-mail:ishigaki.y.aa@m.titech.ac.jp\\
}

\maketitle
\begin{abstract}
We consider the system of equations describing motion of compressible viscoelastic fluids in a whole space. We investigate the large time behavior of solutions around a motionless state, and obtain the $L^p$ decay estimates of solutions for $1<p\leq\infty$, provided that the initial data is sufficiently close to the motionless state. In addition, we clarify the diffusion wave phenomena caused by the sound wave and the elastic shear wave.  
\end{abstract}

\section{Introduction}

This paper studies the compressible viscoelastic system   
\begin{equation}
\left\{\label{system}
\begin{array}{l}
\partial_t\rho+\mathrm{div} (\rho v)=0, \\[1.0ex]
\displaystyle
\rho(\partial_t v+ v\cdot\nabla v)
-\nu{\Delta}v -(\nu+\nu')\nabla \mathrm{div} v
+\nabla  P(\rho)  = \beta^2\mathrm{div}(\rho F{}^\top\! F), \\[1.0ex]
\partial_t F + v\cdot\nabla F =(\nabla v)F
\end{array}
\right.
\end{equation}
in  $\mathbb{R}^3$.
Here $\rho=\rho(x,t)$, 
$v ={}^\top (v^1(x,t), v^2(x,t), v^3(x,t))$, and $F=(F^{jk}(x,t))_{1\leq j,k\leq 3}$ 
 are the unknown density, the velocity field, and the deformation tensor, respectively, at  position $x\in \mathbb{R}^3$ and time $t\geq 0$;
$P=P(\rho) $ is the pressure; $\nu$ and $\nu'$ are the viscosity coefficients satisfying
\[
\nu>0,~2\nu+3\nu'\geq0;
\]
$\beta>0$ is the strength of the elasticity. In particular, if we set $\beta=0$, the system \eqref{system} settles into the usual compressible Navier-Stokes equation. 
We assume that $P'(1)>0$, and we denote $\gamma=\sqrt{P'(1)}$. 

The system \eqref{system} is considered under the initial condition
\begin{align}
(\rho, v, F)|_{t=0}=(\rho_0,v_0,F_0). \label{initialcond1}
\end{align}
We also impose the following conditions
\begin{align}
\left\{
\begin{array}{l}
\mathrm{div}(\rho_0 {}^\top\! F_0)=0,~\rho_0\mathrm{det}F_0=1, \\[1ex]
\displaystyle
\sum_{m=1}^3(F_0^{ml}\partial_{x_m}F_0^{jk}-F_0^{mk}\partial_{x_m}F_0^{jl})=0,
~j,k,l=1,2,3. 
\end{array}
\right.
\label{initialcond2}
\end{align}
According to \cite[Proposition.1]{qianzhang}, these quantities are invariant for $t$:
\begin{align}
\left\{
\begin{array}{l}
\mathrm{div}(\rho{}^\top\! F)=0,~\rho\mathrm{det}F=1,\\
\displaystyle
\sum_{m=1}^3(F^{ml}\partial_{x_m}F^{jk}-F^{mk}\partial_{x_m}F^{jl})=0,~j,k,l=1,2,3,
\end{array}
\right.
~\for~t\geq0. 
\label{constraint1.4}
\end{align}
for $t\geq0$.

The aim of this paper is to study the large time behavior of solutions of the problem \eqref{system}--\eqref{initialcond2} around a motionless state $(\rho,v,F)=(1,0,I)$. Here $I$ is the $3\times3$ identity matrix.

The system \eqref{system} is regarded as one of the basic model describing a motion of viscous compressible fluid with effect of elastic body whose corresponding energy functional is given by
$W(F)=\frac{\beta^2}{2}|F|^2$, called the Hookean linear elasticity. Moreover, we can classify the system \eqref{system} in a quasilinear parabolic-hyperbolic system since the system \eqref{system} is a composite system of the compressible Navier-Stokes equations and a first order hyperbolic system for $F$.  We refer to \cite{gurtin,linliuzhang,sideristhomases} for more physical details. 

In the case $\beta=0$, the large time behavior of the solutions around $(\rho,v)=(1,0)$ has been investigated so far. Matsumura and Nishida \cite{matsumuranishida1979} proved the global existence of the solutions of the problem \eqref{system}--\eqref{initialcond1} provided that the initial perturbation is sufficiently small in $H^3\cap L^1$, and derived the decay estimate:
$$
\|\nabla^k(\phi(t),m(t))\|_{L^2}\leq C(1+t)^{-\frac{3}{4}-\frac{k}{2}}, ~k=0,1,
$$ 
where $(\phi,m)=(\rho-1,\rho v)$.
Hoff and Zumbrun \cite{hoffzumbrun} established the following $L^p~(1\leq p\leq\infty)$ decay estimates in $\mathbb{R}^n,~n\geq 2$:
\begin{gather*}
\|(\phi(t),m(t))\|_{L^p}\leq
\begin{cases}
C(1+t)^{-\frac{n}{2}\left(1-\frac{1}{p}\right)-\frac{n-1}{4}\left(1-\frac{2}{p}\right)}L(t), & 1\leq p < 2, \\
C(1+t)^{-\frac{n}{2}\left(1-\frac{1}{p}\right)}, & 2\leq p \leq\infty,
\end{cases}\\
\end{gather*}
where $L(t)=\log(1+t)$ when $n=2$, and $L(t)=1$ when $n\geq3$. 
Furthermore, the authors of \cite{hoffzumbrun} derived the following asymptotic property:
\begin{gather*}
\left\|\left((\phi(t),m(t))-\left(0,\mathcal{F}^{-1}\left(e^{-\nu|\xi|^2 t}\hat{\mathcal{P}}(\xi)\hat{m}_0\right)\right)\right)\right\|_{L^p}\leq
C(1+t)^{-\frac{n}{2}\left(1-\frac{1}{p}\right)-\frac{n-1}{4}\left(1-\frac{2}{p}\right)}  
\end{gather*}
for $2\leq p \leq \infty$. Here $\hat{\mathcal{P}}(\xi)=I-\frac{{\xi}^\top \xi}{|\xi|^2},~\xi\in\mathbb{R}^n$. According to \cite{kawashimamatsumuranishida}, the solution of the linearized system is expressed as the sum of two terms, one is the incompressible part given by $\mathcal{F}^{-1}\left(e^{-\nu|\xi|^2 t}\hat{\mathcal{P}}(\xi)\hat{m}_0\right)$ which solves the heat equation, and the other is the convolution of the heat kernel and the fundamental solution of the wave equation, called the diffusion wave. The authors of \cite{hoffzumbrun} found that the hyperbolic aspect of the sound wave plays a role of the spreading effect of the wave equation, and the decay rate of the solution becomes slower than the heat kernel when $1\leq p<2$. On the other hand, if $2<p\leq\infty$, the compressible part of the solution $(\phi(t),m(t))-\left(0,\mathcal{F}^{-1}\left(e^{-\nu|\xi|^2 t}\hat{\mathcal{P}}(\xi)\hat{m}_0\right)\right)$ tends to $0$ faster than the heat kernel. See also \cite{kobayashishibata} for the linearized problem.   

We next review the related works in the case $\beta>0$. The local existence of the strong solution of the initial value problem \eqref{system}--\eqref{initialcond2} was shown by Hu and Wang \cite{huwang1}. The global existence of the strong solution of the initial value problem \eqref{system}--\eqref{initialcond2} was proved by Hu and Wang \cite{huwang2}, Qian and Zhang \cite{qianzhang}, and Hu and Wu \cite{huwu}, provided that the initial perturbation $(\rho_0-1,v_0,F_0-I)$ is sufficiently small. Hu and Wu \cite{huwu} also showed that if the initial perturbation $(\rho_0-1,v_0,F_0-I)$ belongs to $L^1(\mathbb{R}^3)\cap H^2(\mathbb{R}^3)$, the $L^p$ decay estimates $\|u(t)\|_{L^p}\leq C(1+t)^{-\frac{3}{2}\left(1-\frac{1}{p}\right)}$ hold for the case $2\leq p\leq6$, via the Fourier splitting method and the Hodge decomposition.
Here $u(t)=(\phi, w, G)=(\rho-1,v,F-I)$.
Li,  Wei and Yao \cite{liweiyao,weiliyao} extended the above result to the case $2\leq p\leq \infty$, and obtained the decay estimates in $L^2$ of higher order derivatives:  
$$
\|\nabla^k(\rho(t)-1,v(t),F(t)-I)\|_{L^2}\leq C(1+t)^{-\frac{3}{4}-\frac{k}{2}},~k=0,1,\ldots,N-1,
$$
provided that $u_0=(\rho_0-1,v_0,F_0-I)$ belongs to $H^N,~N\geq3$, and is small in $L^1\cap H^{3}$. This follows from the diffusive aspect of the system \eqref{system}.
We also refer to \cite{hu,liuwuxuzhang,wugaotan} in recent progresses. 

In view of the results in \cite{hoffzumbrun}, it is expected that the system \eqref{system} has the diffusion wave phenomena affected by the sound wave and the elastic shear wave. In fact, let us consider the linearized system around $(1,0,I)$:
\begin{align}
\partial_t u+Lu=0. \label{LSintro}
\end{align}
Here $L$ is the linearized operator given by
\[
L=
\left(
\begin{array}{@{\ }cc@{\ }cc@{\ } }
0 & \mathrm{div}& 0\\
\gamma^2\nabla & -\nu\Delta-\tilde{\nu}\nabla\mathrm{div} & -\beta^2\mathrm{div} \\
0 & -\nabla & 0 
\end{array}
\right),
\]
where $\tilde\nu=\nu+\nu'$. We then see that the solenoidal part of the velocity $w_s=\mathcal{F}^{-1}(\hat{\mathcal{P}}(\xi)\hat{w})$ satisfies the following linear symmetric parabolic-hyperbolic system:
\[
\left\{
\begin{array}{l}
\partial_t w_s-\nu\Delta w_s-\beta\mathrm{div}\tilde{G}_s=0,\\
\partial_t \tilde{G}_s- \beta\nabla w_s=0,
\end{array}
\right.
\]
where $\tilde{G}_s=\beta\mathcal{F}^{-1}(\hat{\mathcal{P}}(\xi)\hat{G})$,
while the complimentary part $w_c=w-w_s$ solves the following strongly damped wave equation:  
\[
\partial_t^2 w_c-(\beta^2+\gamma^2)\Delta w_c-(\nu+\tilde\nu)\partial_t\Delta w_c=0.
\]
In view of \cite{shibataVW2000}, the solution of linearized system \eqref{LSintro} behaves different to the case $\beta=0$ (\cite{hoffzumbrun,kobayashishibata}) by the additional hyperbolic aspect arising from the elastic shear wave. As a result, the principal part of the linearized system \eqref{LSintro} can be regarded as a system of the strongly damped wave equation.   

In this paper, we shall show that if the initial perturbation $u_0=(\rho_0-1, v_0, F_0-I)$ is sufficiently small in $L^1\cap H^3$, then the global strong solution satisfies the following $L^p$ decay estimate
\[
\|(\rho(t)-1, v(t), F(t)-I)\|_{L^p}\leq C(1+t)^{-\frac{3}{2}\left(1-\frac{1}{p}\right)-\frac{1}{2}\left(1-\frac{2}{p}\right)},~1<p\leq\infty,~t\geq0.
\]  
This result improves the decay rate of the $L^p$ norm of the perturbation $u$ obtained in \cite{huwu,liweiyao} for $p>2$. 

We give an outline of the proof of the main result. Since \eqref{constraint1.4} are nonlinear, straightforward application of the semigroup theory does not work well. To overcome this obstacle, we adopt a material coordinate transform which makes the constraint $\mathrm{div}(\rho{}^\top\! F)=0$ a linear one. We first introduce a displacement vector $\tilde\psi=x-X\in \mathbb{R}^3$ used in \cite{qian,sideristhomases}, where $X=X(x,t)$ is the inverse of the material coordinate. Then we see that $F$ has the form $F-I=\nabla\tilde{\psi}+h(\nabla\tilde{\psi})$. Here $h(\nabla\tilde{\psi})$ is a function satisfying  $h(\nabla\tilde{\psi})=O(|\nabla\tilde\psi|^2), |\nabla\psi|\ll1$.  We next make use of the nonlinear transform $\psi=\tilde{\psi}-(-\Delta)^{-1}\mathrm{div} {}^\top (\phi\nabla\tilde{\psi}+(1+\phi) h(\nabla\tilde{\psi}))$. It turns out that the constraint $\mathrm{div}(\rho {}^\top F)=0$ becomes the linear condition $\phi+\mathrm{tr}(\nabla\psi)=\phi+\mathrm{div}\psi=0$. Furthermore, the decay estimate of the $L^p~(1< p\leq \infty)$ norm of $u=(\phi,w,G)$ is obtained from $U=(\phi,w,\nabla\psi)$. Consequently, the $L^p$ decay estimate can be obtained by employing the following integral equation
\[
U(t)=e^{-tL}U(0)+\int_0^t e^{-(t-s)L}N(U) \mathrm{d}s,
\]
where $N(U)=(N_1(U),N_2(U),N_3(U))$ is a nonlinearity such that $N_1+\mathrm{tr}N_3=0$. We decompose $U$ into the low-frequency part $U_1$ and the high-frequency part $U_\infty$. We then apply the linearized analysis to $U_1$-part, and a variant Matsumura-Nishida energy method \cite{matsumura-nishidaEM} to $U_\infty$-part to establish the result in the case $2\leq p\leq\infty$. On the other hand, for $1<p< 2$, we derive the $L^p$ estimate of $U(t)$ by employing the results in \cite{kobayashishibata, shibataVW2000}. Since the above mentioned nonlinear transformation from $\tilde{\psi}$ to $\psi$ includes the nonlocal operator $(-\Delta)^{-1}$, the case $p=1$ is excluded here. See Remark 4.3 below.

This paper is organized as follows. In Section 2 we introduce some notations and function spaces. In Section 3 we state the main result of this paper on the $L^p$ decay estimates. In Section 4 we reformulate the problem to prove the main result. In Section 5, we give a solution formula of the linearized problem and establish the $L^p$ decay estimates in the case $p\geq2$. In Section 6, we prove the $L^p$ decay estimate in the remaining case $1< p<2$. In the Appendix, we derive the solution formula of the linearized problem. 

\section{Notation}
In this section, we prepare notations and function spaces which will be used throughout the paper. 
$L^p~(1\leq p\leq\infty)$ denotes the usual Lebesgue space on $\mathbb{R}^3$,~and its norm is denoted by $\|\cdot\|_{L^p}$. 
Similarly $W^{m,p} (1\leq p\leq \infty, m\in\{0\}\cup\mathbb{N})$ denotes the $m$-th order $L^p$ Sobolev space on $\mathbb{R}^3$,~and its norm is denoted by $\|\cdot\|_{W^{m,p}}$. We define $H^m=W^{m,2}$ for an integer $m\geq0$. For simplicity, we denote $L^p=L^p\times (L^p)^3\times  (L^p)^9$ (resp. $H^m=H^m\times (H^m)^3\times  (H^m)^9$).

The inner product of $L^2$ is denoted by
$$(f,g):=\int_{\mathbb{R}^3} f(x)\overline{g(x)}dx,~f,g\in L^2.$$
Here the symbol $\overline{\cdot}$ stands for its complex conjugate.
Partial derivatives of a function $u$ in $x_j~(j=1,2,3)$ and $t$ are denoted by $\partial_{x_j}u$ and $\partial_tu$, respectivity. $\Delta$ denotes the usual Laplacian with respect to $x$. For  a multiindex $\alpha=(\alpha_1,\alpha_2,\alpha_3)\in(\{0\}\cup\mathbb{N})^3$ and $\xi={}^\top(\xi_1,\xi_2,\xi_3)\in\mathbb{R}^3$, we denote $\partial_x^\alpha$ and $\xi^\alpha$ by $\partial_x^\alpha =\partial_{x_1}^{\alpha_1} \partial_{x_2}^{\alpha_2} \partial_{x_3}^{\alpha_3}$ and $\xi^\alpha=\xi_1^{\alpha_1}\xi_2^{\alpha_2}\xi_3^{\alpha_3}$, respectivity. For a function $u$ and a nonnegative integer $k$, $\nabla^k u$ stands for $\nabla^k u=\{\partial_x^\alpha u|~|\alpha|=k\}$.

For a scalar valued function $\rho=\rho(x)$, we denote by $\nabla\rho$ its gradient with respect to $x$. For a vector valued function $w=w(x)={}^\top(w^1(x),w^2(x),w^3(x))$, we denote by $\mathrm{div}w$ and $(\nabla w)^{jk}=(\partial_{x_k}w^j)$ its divergence and Jacobian matrix with respect to $x$, respectively. For a $3\times3$-matrix valued function $F=F(x)=(F^{jk}(x))$, we define its divergence $\mathrm{div}F$ and trace $\mathrm{tr}F$ by $(\mathrm{div}F)^j=\sum_{k=1}^3 \partial_{x_k}F^{jk}$ and $\mathrm{tr}F=\sum_{k=1}^3 F^{kk}$, respectively. 

For functions $f=f(x)$ and $g=g(x)$, we denote the convolution of $f$ and $g$ by $f* g$: 
\[
(f* g)(x)=\int_{\mathbb{R}^3} f(x-y)g(y)\mathrm{d}y.
\]

We denote the Fourier transform of a function $f=f(x)$ by $\hat{f}$ or $\mathcal{F}f$:
\[
\hat{f}(\xi)=(\mathcal{F}f)(\xi)=\frac{1}{(2\pi)^{\frac{3}{2}}}\int_{\mathbb{R}^3} f(x)e^{-i\xi\cdot x} \mathrm{d}x~(\xi\in\mathbb{R}^3).
\] 
The Fourier inverse transform is denoted by $\mathcal{F}^{-1}$:
\[
(\mathcal{F}^{-1}f)(x)=\frac{1}{(2\pi)^{\frac{3}{2}}}\int_{\mathbb{R}^3} f(\xi)e^{i\xi\cdot x} \mathrm{d}\xi~(x\in\mathbb{R}^3).
\] 

We recall the Sobolev inequalities.
\begin{lem}\label{Sobolevineq}
The following inequalities hold:
\begin{align*}
\mathrm{(i)}~\|u\|_{L^p}&\leq C\|u\|_{H^1}~\mathrm{for}~2\leq p\leq6,~u\in H^1.\\[1ex]
\mathrm{(ii)}~\|u\|_{L^p}&\leq C\|u\|_{H^2}~\mathrm{for}~2\leq p \leq\infty,~u\in H^2.
\end{align*}
\end{lem}
We will also use the following elementary inequality. (See e.g., \cite[Lemma 3.1]{segal} for the proof.)
\begin{lem}\label{betafunc1}
If $\max\{ a,b\}>1$, then the following estimate holds:
\[
\int_{0}^{t} (1+t-s)^{-a}(1+s)^{-b}\mathrm{d}s\leq C(1+t)^{-\min\{a,b\}},~t\geq0.
\]
\end{lem}

\section{Main Result}
In this section, we state the main result of this paper.

We set
$u(t)=(\phi(t),w(t),G(t))=(\rho(t)-1,v(t),F(t)-I).$
Then $u(t)$ satisfies the following initial value problem
\begin{equation}
\left\{\label{system2}
\begin{array}{l}
\partial_t\phi+\mathrm{div} w=g_1,\\[1ex]
\partial_t  w-\nu{\Delta} w-\tilde \nu\nabla \mathrm{div}  w+\gamma^2\nabla\phi-\beta^2\mathrm{div} G = g_2, \\[1ex]
\partial_t G -\nabla w= g_3,\\[1ex]
\nabla \phi +\mathrm{div}{}^\top\! G = g_4,\\[1ex]
u|_{t=0}=u_0=(\phi_0,w_0,G_0).
\end{array}
\right.
\end{equation}
Here $g_j, j=1, 2, 3, 4,$ denote the nonlinear terms;
\begin{gather*}
\begin{array}{l}
g_1=-\mathrm{div}(\phi w), \\[1ex]
\displaystyle
g_2=-w\cdot\nabla w+\frac{\phi}{1+\phi}(-\nu\Delta w-\tilde{\nu}\nabla\mathrm{div}w+\gamma^2\nabla \phi)-\frac{1}{1+\phi}\nabla Q(\phi)\\[2ex]
\displaystyle
\quad\quad-\frac{\beta^2\phi}{1+\phi}\mathrm{div} G +\frac{\beta^2}{1+\phi}\mathrm{div}(\phi G+G{}^\top\! G+\phi G{}^\top\! G), \\[2ex]
g_3=-w\cdot\nabla G+\nabla w G, \\[1ex]
g_4=-\mathrm{div}(\phi{}^\top\! G),\\[1ex] 
\end{array}   
\end{gather*}
where
\[
Q(\phi)=\phi^2\int_0^1 {P}^{\prime\prime}(1+s\phi)\mathrm{d}s,~\nabla Q=O(\phi)\nabla \phi
\]
for $|\phi|\ll1$. 

We recall the $L^2$ decay estimates obtained in \cite{liweiyao}.

\begin{prop}\label{GE}{\sc $(\cite{liweiyao})$}
Let $u_0\in H^N,~N\geq3$. There is a positive number $\epsilon_0$ such that if $u_0$ satisfies $\|u_0\|_{L^1}+\|u_0\|_{H^3}\leq\epsilon_0$, then there exists a unique solution $u(t)\in C([0,\infty);H^N)$ of the problem $\eqref{system2}$, and $u(t)=(\phi(t),w(t),G(t))$ satisfies 
\begin{gather*}
\| u(t)\|_{H^N}^2+\int_{0}^{t}(\|\nabla\phi(s)\|_{H^{N-1}}^2+\|\nabla w(s)\|_{H^N}^2+\|\nabla G(s)\|_{H^{N-1}}^2)\mathrm{d}s\leq C\|u_0\|_{H^N}^2,
\\
\|\nabla^k u(t)\|_{L^2}\leq C(1+t)^{-\frac{3}{4}-\frac{k}{2}}(\|u_0\|_{L^1}+\|u_0\|_{H^N}) 
\end{gather*}
for $k=0,1,2,\ldots,N-1$ and $t\geq 0$.
\end{prop}

We next state the main result of this paper which reflects an effect of diffusion waves in decay properties.

\begin{thm}\label{mainthm}

$\mathrm{(i)}$ Let $2\leq p\leq\infty$. Assume that $\phi_0,~G_0,$ and $F_0^{-1}$ satisfy $\nabla \phi_0 -\mathrm{div}{}^\top\! (I+G_0)^{-1} =0$ and $F_0^{-1}=\nabla X_0$ for some vector field $X_0$. There is a positive number $\epsilon$ such that if $u_0=(\phi_0,w_0,G_0)$ satisfies $\|u_0\|_{H^3}\leq\epsilon$ and $u_0\in L^1$, then there exists a unique solution $u(t)\in C([0,\infty);H^3)$ of the problem $\eqref{system2}$, and $u(t)=(\phi(t),w(t),G(t))$ satisfies 
\begin{align*}
\|u(t)\|_{L^p}\leq C(p)(1+t)^{-\frac{3}{2}\left(1-\frac{1}{p}\right)-\frac{1}{2}\left(1-\frac{2}{p}\right)}(\|u_0\|_{L^1}+\|u_0\|_{H^3}) 
\end{align*}
uniformly for $t\geq 0$.
Here $C(p)$ is a positive constant depending only on $p$.

$\mathrm{(ii)}$ Let $1<p<2$. Assume that $\phi_0,~G_0,$ and $F_0^{-1}$ satisfy $\nabla \phi_0 -\mathrm{div}{}^\top\! (I+G_0)^{-1}=0$ and $F_0^{-1}=\nabla X_0$ for some vector field $X_0$. There is a positive number $\epsilon_p$ such that if $u_0=(\phi_0,w_0,G_0)$ satisfies $\|u_0\|_{H^3}\leq\epsilon_p$ and $u_0\in L^1$, then there exists a unique solution $u(t)\in C([0,\infty);H^3)$ of the problem $\eqref{system2}$, and $u(t)=(\phi(t),w(t),G(t))$ satisfies 
\begin{align*}
\|u(t)\|_{L^p}\leq C(p)(1+t)^{-\frac{3}{2}\left(1-\frac{1}{p}\right)+\frac{1}{2}\left(\frac{2}{p}-1\right)}(\|u_0\|_{L^1}+\|u_0\|_{L^p}+\|u_0\|_{H^3}) 
\end{align*}
uniformly for $t\geq 0$.
Here $C(p)$ is a positive constant depending only on $p$.
\end{thm}

\begin{rem}
{\rm
Since $\frac{1}{2}\left(1-\frac{2}{p}\right)>0$ for $2< p\leq\infty$, Theorem \ref{mainthm} (i) implies that the $L^p$ norm of the perturbation $u=(\phi, w, G)$ tends to 0 faster than the heat kernel as $t\to\infty$. We thus improve the result in \cite{liweiyao}. Furthermore, due to the elastic force $\beta^2\mathrm{div}(\rho F{}^\top F)$, we discover that if $2< p\leq\infty$, then the decay rate of the $L^p$ norm is faster than the result in \cite{hoffzumbrun} for the  compressible Navier-Stokes equations.
}
\end{rem}

\section{Formulation of the Problem} 
In this section, we rewrite the problem \eqref{system2} into a specific form to prove Theorem \ref{mainthm}.

Let $x=x(X,t)$ be the material coordinate defined by the solution of the flow map:
\begin{numcases}
{}
\frac{dx}{dt}(X,t)=v(x(X,t),t),\nonumber\\
x(X,0)=X,\nonumber
\end{numcases}
and we denote its inverse by $X=X(x,t)$.
According to \cite{gurtin,sideristhomases}, $F$ is defined by $F=\frac{\partial x}{\partial X}$. It is shown in \cite{qian} that its inverse $F^{-1}$ is written as $F^{-1}(x,t)=\nabla X(x,t)$ if $F_0^{-1}$ has the form $F^{-1}_0=\nabla X_0$.  
We set $\tilde{\psi}=x-X$.
Then $\tilde{\psi}$ is a solution of 
\begin{align*}
\partial_t\tilde{\psi}-v=-v\cdot\nabla\tilde{\psi},
\end{align*}
and satisfies
\begin{equation}
G=\nabla\tilde{\psi}+h(\nabla\tilde{\psi}),\label{3.1}
\end{equation}
where $h(\nabla\tilde{\psi})=(I-\nabla\tilde{\psi})^{-1}-I-\nabla\tilde{\psi}$. 

We note that \eqref{3.1} is equivalent to
\begin{equation}
\nabla\tilde{\psi}=I-(I+G)^{-1}.\label{3.2}
\end{equation}

The following estimates hold for $G$ and $\nabla\tilde{\psi}$．

\begin{lem}\label{equivforGpsi}
Assume that $G$ and $\tilde\psi$ satisfy \eqref{3.1}. There is a positive number $\delta_0$ such that if $\|G\|_{H^3}\leq\min\{1,\delta_0\}$, the following inequalities hold:
\begin{align}
&C^{-1}\|\nabla\tilde{\psi}\|_{L^p}\leq\|G\|_{L^p} \leq C\|\nabla\tilde{\psi}\|_{L^p},~1\leq p\leq \infty,  \label{equiv2} \\
&\| \nabla^2\tilde{\psi}\|_{L^2} \leq C\|\nabla G\|_{L^2}, \label{equiv1.1}  \\
&\| \nabla^3\tilde{\psi}\|_{L^2} \leq C(\|\nabla G\|_{H^1}^2+\|\nabla^2 G\|_{L^2}), \label{equiv1.2}  \\
&\| \nabla^4\tilde{\psi}\|_{L^2} \leq C(\|\nabla G\|_{H^1}\|\nabla^2 G\|_{H^1}+\|\nabla^3 G\|_{L^2}). \label{equiv1.3}  
\end{align}
\end{lem}
\textbf{Proof.} If $|G|<1$, \eqref{3.2} implies
\begin{align}
\nabla\tilde\psi=G-\sum_{l=2}^\infty (-G)^l. \label{NExpa}
\end{align}
Let $c>0$ be a positive constant such that $\|G\|_{L^\infty}\leq c\|G\|_{H^2}$.
If $\|G\|_{H^3}\leq\frac{1}{27c}$, then we have $|G|\leq9\|G\|_{L^\infty}\leq9c\|G\|_{H^2}\leq \frac{1}{3}$, and hence
\begin{align*}
\left|\sum_{l=2}^\infty (-G)^l\right|&
\leq  \sum_{l=2}^\infty |G|^{l-1}|G|
\leq \sum_{l=2}^\infty\left(\frac{1}{3}\right)^{l-1}|G|
\leq \frac{1}{2}|G|.
\end{align*}
Therefore, we obtain
\begin{equation}
\left\|\sum_{l=2}^\infty (-G)^l\right\|_{L^p}\leq \frac{1}{2}\|G\|_{L^p}~\mathrm{for}~1\leq p\leq\infty. \label{pfl411}
\end{equation}
Combining \eqref{NExpa} and \eqref{pfl411} yields \eqref{equiv2}.

To prove \eqref{equiv1.1}, we make use of \eqref{3.2}, \eqref{equiv2} and the following formula
\[
\partial_{x_j}( F^{-1})=-F^{-1}\partial_{x_j} F F^{-1},~j=1,2,3.
\]
It then follows that
\begin{align*}
\|\nabla\partial_{x_j} \tilde{\psi}\|_{L^2}&= \|\partial_{x_j} (I+G)^{-1}\|_{L^2} \\
&=\|(I+G)^{-1}\partial_{x_j}G(I+G)^{-1}\|_{L^2}\\
&\leq C\|(I+G)^{-1}\|_{L^\infty}^2\|\partial_{x_j}G\|_{L^2} \\
&\leq C\|\partial_{x_j}G\|_{L^2}.
\end{align*}
This gives \eqref{equiv1.1}.

We next consider \eqref{equiv1.2}.
Since
\begin{align*}
\nabla\partial_{x_j}\partial_{x_k} \tilde{\psi}
&=-\partial_{x_j}\partial_{x_k} (I+G)^{-1}\\
&=-(I+G)^{-1}\partial_{x_j}G(I+G)^{-1}\partial_{x_k}G(I+G)^{-1} \\
&\quad\quad+(I+G)^{-1}\partial_{x_j}\partial_{x_k}G(I+G)^{-1}\\
&\quad\quad\quad-(I+G)^{-1}\partial_{x_k}G(I+G)^{-1}\partial_{x_j}G(I+G)^{-1},
\end{align*}
we have the following estimate by using Lemma 2.1
\begin{align*}
&\|\nabla^3 \tilde{\psi}\|_{L^2} \\
&\leq C(\|(I+G)^{-1}\|_{L^\infty}^3\|\nabla G\|_{L^4}^2+\|(I+G)^{-1}\|_{L^\infty}\|\nabla^2 G\|_{L^2}) \\
&\leq C(\|\nabla G\|_{H^1}^2+\|\nabla^2 G\|_{L^2}).
\end{align*}
We thus obtain \eqref{equiv1.2}.
By a similar computation, we have \eqref{equiv1.3}.
This completes the proof. $\blacksquare$

Based on Lemma \ref{equivforGpsi}, we consider $\tilde{\psi}$ instead of $G$. In terms of $\tilde{U}=(\phi,w,\nabla\tilde{\psi})$, the problem \eqref{system2} is transformed into
\begin{equation}
\left\{\label{system3}
\begin{array}{l}
\partial_t\phi+\mathrm{div} w=f_1,\\
\partial_t  w-\nu{\Delta} w-\tilde \nu\nabla \mathrm{div}  w+\gamma^2\nabla\phi-\beta^2\Delta\tilde{\psi} = f_2, \\
\partial_t \nabla\tilde{\psi} - \nabla w=  f_3,\\
\nabla \phi +\nabla\mathrm{div}\tilde{\psi} = f_4,\\
\tilde{U}|_{t=0}=\tilde{U}_0=(\phi_0,w_0,\nabla\tilde{\psi}_0).
\end{array}
\right.
\end{equation}
Here $f_j, j=1, 2, 3, 4,$ denote the nonlinear terms;
\begin{gather*}
\begin{array}{l}
f_1=g_1, \\[1ex]
\displaystyle
f_2=g_2+\beta^2\mathrm{div}h(\nabla\tilde{\psi}),\\[2ex]
f_3=-\nabla(w\cdot\nabla\tilde{\psi}) , \\[1ex]
f_4=-\mathrm{div}{}^\top\!(\phi\nabla\tilde{\psi}+(1+\phi) h(\nabla\tilde{\psi})).\\[1ex] 
\end{array}   
\end{gather*}

We next introduce $\psi$ by
$\psi=\tilde{\psi}-(-\Delta)^{-1}\mathrm{div} {}^\top (\phi\nabla\tilde{\psi}+(1+\phi) h(\nabla\tilde{\psi}))$, where $(-\Delta)^{-1}=\mathcal{F}^{-1}|\xi|^{-2}\mathcal{F}$, and set $\varPsi=\nabla\psi$. By this transformation, the nonlinear constraint $\nabla\phi+\nabla\mathrm{div}\tilde\psi=f_4$ is transformed into the linear constraint $\phi  +\mathrm{tr}\nabla\psi=0$; and the problem \eqref{system3} is rewritten as
\begin{equation}
\left\{
\begin{array}{l}
\partial_t\phi+\mathrm{div} w=N_1,\\
\partial_t  w-\nu{\Delta} w-\tilde \nu\nabla \mathrm{div}  w+\gamma^2\nabla\phi-\beta^2\mathrm{div} \varPsi = N_2, \\
\partial_t \varPsi-\nabla w= N_3,\\
\phi  +\mathrm{tr}\varPsi=0,~\varPsi=\nabla\psi, \\
U|_{t=0}=U_0=(\phi_0,w_0,\varPsi_0).
\end{array}
\right.
\label{problem4}
\end{equation}
Here $N_j, j=1, 2, 3,$ denote the nonlinear terms;
\begin{gather*}
\begin{array}{l}
N_1=f_1, \\[1ex]
\displaystyle
N_2=f_2-\beta^2\mathrm{div}{}^\top (\phi\nabla\tilde{\psi}+(1+\phi) h(\nabla\tilde{\psi})), \\[2ex]
N_3=-\nabla(w\cdot\nabla\tilde{\psi})-\nabla(-\Delta)^{-1}\nabla\mathrm{div}(\phi w)-\nabla(-\Delta)^{-1}\nabla\mathrm{div}(w\cdot\nabla\tilde{\psi}). \\[1ex]
\end{array}   
\end{gather*}
We note that $N_1$ and $N_3$ satisfy $N_1+\mathrm{tr}N_3=0$.
The relations between $\psi$ and $\tilde{\psi}$ are given as follows.

\begin{lem}\label{psitildepsi}

$\mathrm{(i)}$ Let $\tilde{U}_0$ and $U_0$  be the ones as in \eqref{system3} and \eqref{problem4}, respectively. If $\phi_0$ and $\tilde{\psi}_0$ satisfy $\nabla \phi_0 +\nabla\mathrm{div}\tilde{\psi}_0=0$, then it holds $U_0=\tilde{U}_0=(\phi_0,w_0,\nabla\tilde{\psi}_0)$.

$\mathrm{(ii)}$ There is a positive number $\delta_0$ such that the following assertion holds true. Let 
\[
\phi\in  C([0,\infty);H^{3}),~\psi\in   C([0,\infty);H^{4}).
\]
If $\|\phi\|_{C([0,\infty);H^3)}+\|\psi\|_{C([0,\infty);H^4)}\leq\delta_0$, then there uniquely exists $\tilde\psi\in C([0,\infty);H^{4})$ such that
\begin{align}
&\|\tilde\psi\|_{C([0,\infty);H^{4})}\leq\sqrt{\delta_0}, \nonumber\\
&\tilde\psi=\psi+(-\Delta)^{-1}\mathrm{div} {}^\top (\phi\nabla\tilde{\psi}+(1+\phi) h(\nabla\tilde{\psi})). \label{Lem4.2eq1}
\end{align}
$\mathrm{(iii)}$ Let $1< p <\infty$. There is a positive number $\delta_p$ such that if $\|\phi\|_{C([0,\infty);H^3)}+\|\nabla\tilde\psi\|_{C([0,\infty);H^3)}\leq\min\{\delta_0,\delta_p\}$, the following inequalities hold for $t\geq0$:
\begin{align}
&C_p^{-1}\|\nabla\tilde{\psi}(t)\|_{L^p}\leq\|\nabla\psi(t)\|_{L^p} \leq C_p\|\nabla\tilde{\psi}(t)\|_{L^p}. \label{equiv3} 
\end{align}
$\mathrm{(iv)}$ There is a positive number $\delta_1$ such that if $\|\phi\|_{C([0,\infty);H^3)}+\|\nabla\tilde\psi\|_{C([0,\infty);H^3)}\leq\min\{\delta_0,\delta_1\}$, the following inequalities hold for $t\geq0$:
\begin{align}
&
\begin{array}{l}
\|\nabla\tilde{\psi}(t)\|_{L^\infty}\leq C(\|\phi(t)\|_{L^\infty} +\|\nabla\psi(t)\|_{L^\infty} )\\[1ex]
\quad\quad\quad\quad\quad\quad\quad
+C(\|\nabla \phi(t)\|_{H^1}+\|\nabla^2\tilde{\psi}(t)\|_{H^1})^2,
\end{array}
\label{equiv4}\\
&
\begin{array}{l}
\|\nabla^2\psi(t)\|_{L^2}\leq C\|\nabla^2\tilde{\psi}(t)\|_{L^2}\\[1ex]
\quad\quad\quad\quad\quad\quad\quad
+C(\|\phi(t)\|_{H^2}+\|\nabla\tilde{\psi}(t)\|_{H^2})\|\nabla\tilde{\psi}(t)\|_{H^2}, 
\end{array}
\label{equiv5.1}\\
&
\begin{array}{l}
\|\nabla^3\psi(t)\|_{L^2}\leq C(1+\|\phi(t)\|_{H^2}+\|\nabla\tilde{\psi}(t)\|_{H^2})\|\nabla^3\tilde{\psi}(t)\|_{L^2}\\[1ex]
\quad\quad\quad\quad\quad\quad\quad
+C(\|\nabla\phi(t)\|_{H^1}+\|\nabla^2\tilde{\psi}(t)\|_{H^1})\|\nabla\tilde{\psi}(t)\|_{H^2}, 
\end{array}
\label{equiv5.2}
\\
&\|\nabla^4\psi(t)\|_{L^2}\leq C\|\nabla^4\tilde{\psi}(t)\|_{L^2}+C(\|\phi(t)\|_{H^3}+\|\nabla\tilde{\psi}(t)\|_{H^3})^2. 
\label{equiv5.3}
\end{align}
\end{lem}
\textbf{Proof.}
(i) The condition $\nabla\phi_0+\nabla\mathrm{div}\tilde{\psi}_0=0$ leads to $\phi_0\nabla\tilde{\psi}_0+(1+\phi_0) h(\nabla\tilde{\psi}_0)=0$. Therefore, we have $U_0=\tilde{U}_0=(\phi_0,w_0,\nabla\tilde{\psi}_0)$.

(ii) We set $\Gamma(\tilde\psi)=\psi+(-\Delta)^{-1}\mathrm{div} {}^\top (\phi\nabla\tilde{\psi}+(1+\phi) h(\nabla\tilde{\psi}))$ and $\mathcal{B}_{\sqrt{\delta_0}}=\{f\in C([0,\infty),H^4)|~\|f\|_{C([0,\infty),H^4)}\leq \sqrt{\delta_0}\}$. We then see that if $C_1\sqrt{\delta_0}\leq1$, then $\Gamma$ is a mapping of $\mathcal{B}_{\sqrt{\delta_0}}$ into $\mathcal{B}_{\sqrt{\delta_0}}$. Indeed, since
\begin{align*}
h(\nabla\tilde{\psi})=&((I-\nabla\tilde{\psi})^{-1}-I)\nabla\tilde{\psi},
\end{align*}
we see that if $|\nabla\tilde{\psi}|<1$, then
\begin{align*}
h(\nabla\tilde{\psi})=&\sum_{m=2}^\infty (\nabla\tilde{\psi})^m,
\end{align*}
and hence,
\begin{align*}
(I-\nabla\tilde{\psi})^{-1}-I=&\sum_{m=1}^\infty (\nabla\tilde{\psi})^m.
\end{align*}
Furthermore, since
\begin{align*}
\partial_{x_j} ((I-\nabla\tilde{\psi})^{-1})=&(I-\nabla\tilde{\psi})^{-1}\nabla\partial_{x_j}\tilde{\psi}(I-\nabla\tilde{\psi})^{-1},
\end{align*}
we have
\begin{align*}
\|(-\Delta)^{-1}\mathrm{div}{}^\top h(\nabla\tilde{\psi})\|_{L^2}
&\leq C(\|h(\nabla\tilde{\psi})\|_{L^1}+\|h(\nabla\tilde{\psi})\|_{L^2}) \\
&\leq C\sum_{m=2}^\infty (\|(\nabla\tilde{\psi})^m\|_{L^1}+\|(\nabla\tilde{\psi})^m\|_{L^2}) \\
&\leq C\|\nabla\tilde{\psi}\|_{H^2}^2,
\end{align*}
and similarly,
\begin{align*}
\|\nabla(-\Delta)^{-1}\mathrm{div}{}^\top h(\nabla\tilde{\psi})\|_{L^2}
&\leq C\|\nabla\tilde{\psi}\|_{H^2}\|\nabla\tilde{\psi}\|_{L^2}.
\end{align*}
As for the estimate of the second order derivative of $(-\Delta)^{-1}\mathrm{div}{}^\top h(\nabla\tilde\psi)$, since
\begin{align*}
\partial_{x_j}\partial_{x_k} ((I-\nabla\tilde{\psi})^{-1})
=&(I-\nabla\tilde{\psi})^{-1}\nabla\partial_{x_j}\tilde{\psi}(I-\nabla\tilde{\psi})^{-1}\nabla\partial_{x_k}\tilde{\psi}(I-\nabla\tilde{\psi})^{-1} \\
&+(I-\nabla\tilde{\psi})^{-1}\nabla\partial_{x_j}\partial_{x_k}\tilde{\psi}(I-\nabla\tilde{\psi})^{-1}\\
&+(I-\nabla\tilde{\psi})^{-1}\nabla\partial_{x_k}\tilde{\psi}(I-\nabla\tilde{\psi})^{-1}\nabla\partial_{x_j}\tilde{\psi}(I-\nabla\tilde{\psi})^{-1},
\end{align*}
we have
\begin{align*}
&\|\nabla^2(-\Delta)^{-1}\mathrm{div}{}^\top h(\nabla\tilde{\psi})\|_{L^2}\\
&\leq C\|\nabla h(\nabla\tilde{\psi})\|_{L^2} \\
&\leq C(\|(I-\nabla\tilde{\psi})^{-1}-I\|_{L^\infty}\|\nabla^2\tilde{\psi}\|_{L^2}+\|(I-\nabla\tilde{\psi})^{-1}\|_{L^\infty}^2\|\nabla\tilde{\psi}\|_{L^2}\|\nabla^2\tilde{\psi}\|_{L^2}) \\
&\leq C\|\nabla\tilde{\psi}\|_{H^2}\|\nabla^2\tilde{\psi}\|_{L^2}.
\end{align*}
Similarly, one can show that
\begin{align*}
\|\nabla^3(-\Delta)^{-1}\mathrm{div}{}^\top h(\nabla\tilde{\psi})\|_{L^2} 
&\leq C(\|\nabla\tilde{\psi}\|_{H^2}\|\nabla^3\tilde{\psi}\|_{L^2}+\|\nabla^2\tilde{\psi}\|_{H^1}^2),
\end{align*}
and
\[
\|\nabla^4(-\Delta)^{-1}\mathrm{div}{}^\top h(\nabla\tilde{\psi})\|_{L^2}\leq C\|\nabla\tilde{\psi}\|_{H^3}^2.
\]
It then follows that if $\tilde\psi\in\mathcal{B}_{\sqrt{\delta_0}}$, then
\begin{align*}
\begin{array}{l}
\|\Gamma(\tilde\psi)\|_{H^4}\\[1ex]
\leq C(\|\psi\|_{H^4}+\|(-\Delta)^{-1}\mathrm{div} {}^\top (\phi\nabla\tilde{\psi}+(1+\phi) h(\nabla\tilde{\psi}))\|_{H^4})\\[1ex]
\leq C(\|\psi\|_{H^4}+\| \phi\nabla\tilde{\psi}+(1+\phi) h(\nabla\tilde{\psi})\|_{L^1}) \\[1ex]
\quad\quad
+C\| \phi\nabla\tilde{\psi}+(1+\phi) h(\nabla\tilde{\psi})\|_{H^3}\\[1ex]
\leq C(\|\psi\|_{H^4}+\| \phi\|_{H^3}\|\nabla\tilde{\psi}\|_{H^3}+(1+\|\phi\|_{H^3})\|\nabla\tilde{\psi}\|_{H^3}^2) \\[1ex]
\leq C_1\delta_0 \\[1ex]
\leq \sqrt{\delta_0}.
\end{array}
\end{align*}
Therefore, $\Gamma(\tilde\psi)$ belongs to $\mathcal{B}_{\sqrt{\delta_0}}$.

We next claim that if $\tilde\psi_j\in\mathcal{B}_{\sqrt{\delta_0}}$ $(j=1,2)$, then
\begin{align}
\|\Gamma(\tilde\psi_{1})-\Gamma(\tilde\psi_{2})\|_{C([0,\infty);H^4)} \leq C_2\sqrt{\delta}_0\|\tilde\psi_{1}-\tilde\psi_2\|_{C([0,\infty);H^4)}. \label{Lem42ii-ineq2}
\end{align}
To show this, we first have
\begin{align*}
&\|\Gamma(\tilde\psi_1)-\Gamma(\tilde\psi_{2})\|_{H^4} \\
&\leq C(\|\phi\nabla(\tilde{\psi}_{1}-\tilde{\psi}_{2})\|_{H^3} +\|(1+\phi) (h(\nabla\tilde{\psi}_{1})-h(\nabla\tilde{\psi}_2))\|_{H^3} \\
&\leq C\delta_0\|\tilde\psi_{1}-\tilde\psi_{2}\|_{H^4}+C(1+\delta_0)\| h(\nabla\tilde{\psi}_{1})-h(\nabla\tilde{\psi}_2)\|_{H^3}.
\end{align*}
As for the second term on the right-hand side, since
\[
\begin{array}{l}
h(\nabla\tilde{\psi}_{1})-h(\nabla\tilde{\psi}_2)\\[1ex]
=((I-\nabla\tilde{\psi}_{1})^{-1}-I)\nabla(\tilde{\psi}_{1}-\tilde{\psi}_{2})\\[1ex]
\quad
+\nabla(\tilde{\psi}_{1}-\tilde{\psi}_{2})((I-\nabla\tilde{\psi}_{2})^{-1}-I)\\[1ex]
\quad\quad
+((I-\nabla\tilde{\psi}_{1})^{-1}-I)\nabla(\tilde{\psi}_{1}-\tilde{\psi}_{2})((I-\nabla\tilde{\psi}_{2})^{-1}-I),\\[1ex]
\end{array}
\]
we see that
\[
\| h(\nabla\tilde{\psi}_{1})-h(\nabla\tilde{\psi}_2)\|_{H^3}\leq C(\sqrt{\delta}_0+\delta_0)\|\tilde\psi_{1}-\tilde\psi_{2}\|_{H^4}.
\]
Therefore, we arrive at \eqref{Lem42ii-ineq2}.
Taking $\delta_0$ small such that $C_2\sqrt{\delta}_0<1$, we conclude that $\Gamma$ is a contraction map in $\mathcal{B}_{\sqrt{\delta_0}}$. By the contraction mapping principle, we observe that there exists a unique $\tilde\psi\in \mathcal{B}_{\sqrt{\delta_0}}$ such that $\tilde\psi=\Gamma(\tilde\psi)$. This indicates the unique existence of $\tilde\psi$ satisfying \eqref{Lem4.2eq1}.

(iii) We assume that $\|\phi\|_{C([0,\infty);H^3)}+\|\tilde\psi\|_{C([0,\infty);H^4)}\leq \delta$ with some small number $0<\delta<1$ to be determined later. Since the Riesz operator $\mathcal{R}_jf=\mathcal{F}^{-1}\left[\frac{\xi_j}{|\xi|}\hat{f}\right]$ is bounded from $L^p$ to $L^p$ for $1<p<\infty$, we have
\begin{align*}
&\|\nabla(-\Delta)^{-1}\mathrm{div} {}^\top (\phi\nabla\tilde{\psi}+(1+\phi) h(\nabla\tilde{\psi}))\|_{L^p}\\
&=\left\|\mathcal{F}^{-1}\left[\mathcal{F}(^\top(\phi\nabla\tilde{\psi}+(1+\phi) h(\nabla\tilde{\psi})))\frac{\xi^\top\xi}{|\xi|^2}\right]\right\|_{L^p} \\
&\leq C_p\|\phi\nabla\tilde{\psi}+(1+\phi) h(\nabla\tilde{\psi})\|_{L^p}  \\
&\leq C_p(\|\phi\|_{H^2}\|\nabla\tilde{\psi}\|_{L^p}+\| h(\nabla\tilde{\psi})\|_{L^p}+\|\phi\|_{H^2}\| h(\nabla\tilde{\psi})\|_{L^p}).
\end{align*}
We see from $h(\nabla\tilde\psi)=\sum_{m=2}^\infty (\nabla\tilde\psi)^m$ that
\[
\| h(\nabla\tilde{\psi})\|_{L^p}\leq C\|\nabla\tilde{\psi}\|_{L^\infty}\|\nabla\tilde{\psi}\|_{L^p}.
\]
This leads to the estimate
\[
\|\nabla(-\Delta)^{-1}\mathrm{div} {}^\top (\phi\nabla\tilde{\psi}+(1+\phi) h(\nabla\tilde{\psi}))\|_{L^p}\leq C_p(\delta+\delta^2)\|\nabla\tilde{\psi}\|_{L^p}.
\]
By taking $\delta$ small such that $C_p(\delta+\delta^2)\leq\frac{1}{2}$, we obtain \eqref{equiv3}.

(iv) We assume that $\|\phi\|_{C([0,\infty);H^3)}+\|\tilde\psi\|_{C([0,\infty);H^4)}\leq \delta$ with some small number $0<\delta<1$ to be determined later. It follows from the Sobolev inequality and the Plancherel theorem that
\begin{align*}
&\|\nabla(-\Delta)^{-1}\mathrm{div} {}^\top (\phi\nabla\tilde{\psi}+(1+\phi) h(\nabla\tilde{\psi}))\|_{L^\infty}\\
&=\left\|\mathcal{F}^{-1}\left[\mathcal{F}(^\top(\phi\nabla\tilde{\psi}+(1+\phi) h(\nabla\tilde{\psi})))\frac{\xi^\top\xi}{|\xi|^2}\right]\right\|_{L^\infty} \\
&\leq C\left\|\mathcal{F}^{-1}\left[\mathcal{F}(^\top(\phi\nabla\tilde{\psi}+(1+\phi) h(\nabla\tilde{\psi})))\frac{\xi^\top\xi}{|\xi|^2}\right]\right\|_{H^2} \\
&\leq C(\|\phi\nabla\tilde{\psi}\|_{H^2}+\|h(\nabla\tilde{\psi})\|_{H^2}+\|\phi h(\nabla\tilde{\psi})\|_{H^2}).
\end{align*}
Since
\begin{align*}
\|\phi\nabla\tilde{\psi}\|_{H^2}&\leq C(\|\phi\|_{H^2}+\|\nabla\tilde{\psi}\|_{H^2})(\|\phi\|_{L^\infty}+\|\nabla\tilde{\psi}\|_{L^\infty})\\
&\quad+C(\|\nabla \phi\|_{H^1}+\|\nabla^2 \tilde{\psi}\|_{H^1})^2,  \\
\|h(\nabla\tilde{\psi})\|_{H^2}&\leq C\|\nabla\tilde{\psi}\|_{H^2}\|\nabla\tilde{\psi}\|_{L^\infty}+C\|\nabla^2\tilde{\psi}\|_{H^1}^2,  \\
\|\phi h(\nabla\tilde{\psi})\|_{H^2}&\leq C\|\phi\|_{H^2}(\|\nabla\tilde{\psi}\|_{H^2}\|\nabla\tilde{\psi}\|_{L^\infty}+\|\nabla^2\tilde{\psi}\|_{H^1}^2), 
\end{align*}
we have
\begin{align*}
\|\nabla\tilde{\psi}\|_{L^\infty}
&=\|\nabla\psi+\nabla(-\Delta)^{-1}\mathrm{div} {}^\top (\phi\nabla\tilde{\psi}+(1+\phi) h(\nabla\tilde{\psi}))\|_{L^\infty} \\
&\leq C(\delta+\delta^2)\|\nabla\tilde{\psi}\|_{L^\infty}\\
&\quad+C(\|\phi\|_{L^\infty}+\|\nabla\psi\|_{L^\infty})+C(\|\nabla \phi\|_{H^1}+\|\nabla^2 \tilde{\psi}\|_{H^1})^2.
\end{align*}
By taking $\delta$ small such that $C(\delta+\delta^2)\leq\frac{1}{2}$, we arrive at \eqref{equiv4}. We can derive \eqref{equiv5.1}--\eqref{equiv5.3} as in the proof of (ii). This completes the proof. $\blacksquare$

\begin{rem}
{\rm
Due to the restriction $p>1$ in Lemma \ref{psitildepsi} (iii), the decay estimate of $L^1$ norm of $u(t)$ is excluded in Theorem \ref{mainthm}. 
}
\end{rem}

\section{Proof of Theorem \ref{mainthm} (i)} 
In this section, we prove Theorem \ref{mainthm} (i). 
The global existence and the $L^2$ decay estimates of higher order derivatives are guaranteed by Proposition \ref{GE}. Hence we focus on the derivation of the $L^p$ decay estimates except the case $p=2$.
In view of Lemmata \ref{equivforGpsi}--\ref{psitildepsi} and the interpolation inequality:
$
\|u(t)\|_{L^p}\leq\|u(t)\|_{L^2}^{\frac{2}{p}}\|u(t)\|_{L^\infty}^{1-\frac{2}{p}}~(2\leq p\leq\infty),
$
it suffices to obtain the $L^\infty$ decay estimate of $U=(\phi,w,\varPsi)$.

The problem \eqref{problem4} is written in the form:
\begin{equation}
\left\{\label{problem5}
\begin{array}{l}
\partial_t U+LU=N, \\
\phi +\mathrm{div}\psi = 0,\\
U|_{t=0}=U_0,
\end{array}
\right.
\end{equation}
where
$$
L=
\left(
\begin{array}{@{\ }cc@{\ }cc@{\ } }
0 & \mathrm{div}& 0\\
\gamma^2\nabla & -\nu\Delta-\tilde{\nu}\nabla\mathrm{div} & -\beta^2\mathrm{div} \\
0 & -\nabla & 0 
\end{array}
\right),~
N=\left(
\begin{array}{l}
N_1 \\
N_2 \\
N_3
\end{array}
\right).
$$

Theorem \ref{mainthm} (i) is proved by combining Lemma \ref{equivforGpsi}, Lemma \ref{psitildepsi} and the following $L^\infty$ decay estimate of $U(t)$. 
\begin{prop}\label{LinftyestimateofU}
There exists a positive number $\delta_0$ such that if $\|u_0\|_{L^1}+\|u_0\|_{H^3}\leq\delta_0$, then the following inequality 
$$
\|U(t)\|_{L^\infty}\leq C(1+t)^{-2}(\|u_0\|_{L^1}+\|u_0\|_{H^3})
$$
holds for $t\geq 0$.
\end{prop}

To prove Proposition \ref{LinftyestimateofU}, we first give the following $L^2$ decay estimates for $\nabla^k U(t)$.
\begin{prop}\label{H2estimate}
There exists a positive number $\delta_0$ such that if $\|u_0\|_{L^1}+\|u_0\|_{H^3}\leq\delta_0$, then the following inequality 
$$
\|\nabla^k U(t)\|_{L^2}\leq C(1+t)^{-\frac{3}{4}-\frac{k}{2}}(\|u_0\|_{L^1}+\|u_0\|_{H^3})
$$
holds for $k=0,1,2$ and $t\geq 0$.
\end{prop}
Proposition \ref{H2estimate} follows from Proposition \ref{GE}, Lemma \ref{equivforGpsi} and Lemma \ref{psitildepsi}.

We next investigate the linearized problem 
\begin{equation}
\left\{\label{linearizedproblem5}
\begin{array}{l}
\partial_t U+LU=0, \\
\phi +\mathrm{div}\psi = 0,\\
U|_{t=0}=U_0.
\end{array}
\right.
\end{equation}
We denote by $e^{-tL}$ the semigroup generated by $-L$. The solution of \eqref{linearizedproblem5} is written as $U(t)=e^{-tL}U_0$.

To investigate the large time behavior of $U(t)=e^{-tL}U_0$, we take the Fourier transform with respect to $x$. We then obtain
\begin{equation}
\left\{\label{Flinearizedproblem5}
\begin{array}{l}
\partial_t \hat{U}+\hat{L}_\xi\hat{U}=0, \\
\hat\phi +i\xi\cdot\hat\psi = 0,\\
\hat{U}|_{t=0}=\hat{U}_0,
\end{array}
\right.
\end{equation}
where
$$
\hat{L}_\xi\hat{U}=
\left(
\begin{array}{@{\ }cc@{\ }cc@{\ } }
i\xi\cdot\hat{w}\\
i\gamma^2\hat\phi\xi+(\nu|\xi|^2I+\tilde\nu\xi{}^\top\xi)\hat{w}-i\beta^2\hat\varPsi\xi \\
-i\hat{w}{}^\top\xi
\end{array}
\right).
$$

We have the following expression of $e^{-t\hat{L}_\xi}\hat{U}_0$.
\begin{lem}\label{solform}
If $|\xi|\neq0,\frac{\beta}{\nu},\frac{\sqrt{\beta^2+\gamma^2}}{\nu+\tilde\nu}$, the solution of \eqref{Flinearizedproblem5} is written as 
\begin{align}
\left(
\begin{array}{l}
\hat\phi(\xi,t) \\
\hat{w}(\xi,t) \\
\hat{\varPsi}(\xi,t)
\end{array}
\right)=
\left(
\begin{array}{@{\ }cc@{\ }cc@{\ } }
\hat{K}^{11}(\xi,t) & \hat{K}^{12}(\xi,t)& \hat{K}^{13}(\xi,t)\\
\hat{K}^{21}(\xi,t) & \hat{K}^{22}(\xi,t)& \hat{K}^{23}(\xi,t)\\
\hat{K}^{31}(\xi,t) & \hat{K}^{32}(\xi,t)& \hat{K}^{33}(\xi,t)
\end{array}
\right)
\left(
\begin{array}{l}
\hat\phi_0(\xi) \\
\hat{w}_0(\xi) \\
\hat\varPsi_0(\xi)
\end{array}
\right).
\label{solformula}
\end{align}
Here
\begin{align*}
\hat{K}^{11}(\xi,t) =&\frac{\mu_3(\xi)e^{\mu_4(\xi)t}-\mu_4(\xi)e^{\mu_3(\xi)t}}{\mu_3(\xi)-\mu_4(\xi)},\\
\hat{K}^{12}(\xi,t)=&-i\frac{e^{\mu_3(\xi)t}-e^{\mu_4(\xi)t}}{\mu_3(\xi)-\mu_4(\xi)}{}^\top\xi, \\
\hat{K}^{13}(\xi,t)=&0, 
\end{align*}
\begin{align*}
\hat{K}^{21}(\xi,t)=&-i\gamma^2\frac{e^{\mu_3(\xi)t}-e^{\mu_4(\xi)t}}{\mu_3(\xi)-\mu_4(\xi)}\xi, \\
\hat{K}^{22}(\xi,t) =&\frac{\mu_1(\xi)e^{\mu_1(\xi)t}-\mu_2(\xi)e^{\mu_2(\xi)t}}{\mu_1(\xi)-\mu_2(\xi)}\left(I-\frac{{\xi}^\top \xi}{|\xi|^2}\right) \\
&\quad+\frac{\mu_3(\xi)e^{\mu_3(\xi)t}-\mu_4(\xi)e^{\mu_4(\xi)t}}{\mu_3(\xi)-\mu_4(\xi)}\frac{{\xi}^\top \xi}{|\xi|^2},\\
\hat{K}^{31}(\xi,t)=&0,\\
\hat{K}^{33}(\xi,t)=&\frac{\mu_1(\xi)e^{\mu_2(\xi)t}-\mu_2(\xi)e^{\mu_1(\xi)t}}{\mu_1(\xi)-\mu_2(\xi)}\left(I-\frac{{\xi}^\top \xi}{|\xi|^2}\right)\\
&\quad+\frac{\mu_3(\xi)e^{\mu_4(\xi)t}-\mu_4(\xi)e^{\mu_3(\xi)t}}{\mu_3(\xi)-\mu_4(\xi)}\frac{{\xi}^\top \xi}{|\xi|^2}; 
\end{align*}
$\hat{K}^{23}(\xi,t)\hat{\varPsi}_0(\xi)$ and $\hat{K}^{32}(\xi,t)\hat{w}_0(\xi)$ are defined by
\begin{align*}
\hat{K}^{23}(\xi,t)\hat{\varPsi}_0(\xi)=&i\beta^2\frac{e^{\mu_1(\xi)t}-e^{\mu_2(\xi)t}}{\mu_1(\xi)-\mu_2(\xi)}\left(I-\frac{{\xi}^\top \xi}{|\xi|^2}\right)\hat{\varPsi}_0(\xi)\xi \\
&\quad+i\beta^2\frac{e^{\mu_3(\xi)t}-e^{\mu_4(\xi)t}}{\mu_3(\xi)-\mu_4(\xi)}\frac{{\xi}^\top \xi}{|\xi|^2}\hat{\varPsi}_0(\xi)\xi, 
\end{align*}
\begin{align*}
\hat{K}^{32}(\xi,t)\hat{w}_0(\xi)=&i\frac{e^{\mu_1(\xi)t}-e^{\mu_2(\xi)t}}{\mu_1(\xi)-\mu_2(\xi)}\left(I-\frac{{\xi}^\top \xi}{|\xi|^2}\right)\hat{w}_0(\xi){}^\top\xi \\
&\quad+i\frac{e^{\mu_3(\xi)t}-e^{\mu_4(\xi)t}}{\mu_3(\xi)-\mu_4(\xi)}\frac{{\xi}^\top \xi}{|\xi|^2}\hat{w}_0(\xi){}^\top\xi,
\end{align*}
where $\mu_j(\xi),~ j=1,2,3,4$, are given by
\begin{align*}
\mu_1(\xi)&=\frac{-\nu|\xi|^2+\sqrt{\nu^2|\xi|^4-4\beta^2|\xi|^2} }{2},\\
\mu_2(\xi)&=\frac{-\nu|\xi|^2-\sqrt{\nu^2|\xi|^4-4\beta^2|\xi|^2} }{2},\\
\mu_3(\xi)&=\frac{-(\nu+\tilde{\nu})|\xi|^2+\sqrt{(\nu+\tilde{\nu})^2|\xi|^4-4(\beta^2+\gamma^2)|\xi|^2} }{2},\\
\mu_4(\xi)&=\frac{-(\nu+\tilde{\nu})|\xi|^2-\sqrt{(\nu+\tilde{\nu})^2|\xi|^4-4(\beta^2+\gamma^2)|\xi|^2} }{2}.
\end{align*} 
\end{lem}
The proof of Lemma \ref{solform} will be given in Appendix.

The solution $U(t)=e^{-tL}U_0$ is thus given by
\begin{align*}
U(t)=e^{-tL}U_0=\mathcal{F}^{-1}e^{-t\hat{L}_\xi}\hat{U_0}. 
\end{align*}

To study the asymptotic behavior of $U(t)$, we will make use of the following properties of $\mu_j~(j=1,2,3,4)$:
\begin{gather*}
\mu_j(\xi)^2+\nu|\xi|^2\mu_j(\xi)+\beta^2|\xi|^2=0,~j=1,2,\\
\mu_j(\xi)\sim -\frac{\nu}{2}|\xi|^2+i(-1)^{j+1}\beta|\xi|,~\mathrm{for}~|\xi|\ll1,~j=1,2, \\
\mu_1(\xi)\sim-\frac{\beta}{\nu},~\mu_2(\xi)\sim-\nu|\xi|^2,~\mathrm{for}~|\xi|\gg1,\\
\mu_j(\xi)^2+(\nu+\tilde{\nu})|\xi|^2\mu_j(\xi)+(\beta^2+\gamma^2)|\xi|^2=0,~j=3,4,\\
\mu_j(\xi)\sim-\frac{\nu+\tilde{\nu}}{2}|\xi|^2+i(-1)^{j+1}\sqrt{\beta^2+\gamma^2}|\xi|,~\mathrm{for}~|\xi|\ll1,~j=3,4, \\
\mu_3(\xi)\sim-\frac{\beta^2+\gamma^2}{\nu+\tilde{\nu}},~\mu_4(\xi)\sim-(\nu+\tilde{\nu})|\xi|^2,~\mathrm{for}~|\xi|\gg1. 
\end{gather*}

We decompose the solution $U(t)$ of the problem \eqref{problem5} into its low and high frequency parts. Let $\hat\varphi_1, \hat\varphi_\infty\in C^\infty(\mathbb{R}^3)$ be  cut-off functions such that
\[
\hat\varphi_1(\xi)=
\begin{cases}
1 & |\xi|\leq \frac{M_1}{2}, \\
0 & |\xi|\geq \frac{M_1}{\sqrt{2}},
\end{cases}
~\hat\varphi_1(-\xi)=\hat\varphi_1(\xi),
\]
\[
\hat\varphi_\infty(\xi)=1-\hat\varphi_1(\xi),
\]
where
\[
M_1=\min\left\{\frac{\beta}{\nu},\frac{\sqrt{\beta^2+\gamma^2}}{\nu+\tilde\nu}\right\}.
\]
We define the operators $P_1$ and $P_\infty$ on $L^2$ by
\[
P_1u=\mathcal{F}^{-1}(\hat\varphi_1\hat{u}),~P_\infty u=\mathcal{F}^{-1}(\hat\varphi_\infty\hat{u})~\mathrm{for}~u\in L^2.
\]

\begin{lem}\label{lemP1Pinfty}
$P_j~(j=1,\infty)$ have the following properties.

$(i)$ $P_1+P_\infty=I$.

$(ii)$ $\partial_x^\alpha P_1=P_1\partial_x^\alpha$, $\|\partial_x^\alpha P_1 f\|_{L^2}\leq C_\alpha \|f\|_{L^2}$ for $\alpha\in(\{0\}\cup\mathbb{N})^3$ and $f\in L^2$.

$(iii)$ $\partial_x^\alpha P_\infty=P_\infty\partial_x^\alpha$, $\|\partial_x^\alpha P_\infty f\|_{L^2}\leq C \|\nabla\partial_x^\alpha P_\infty f\|_{L^2}$ for $\alpha\in(\{0\}\cup\mathbb{N})^3$ with $|\alpha|=k\geq0$ and $f\in H^{k+1}$.
\end{lem}
Lemma \ref{lemP1Pinfty} immediately follows from the definitions of $P_j,~j=1,\infty,$ and the Plancherel theorem. We omit the proof. 
 
The solution $U(t)$ of \eqref{problem5} is decomposed as
\[
U(t)=U_1(t)+U_\infty(t),~U_1(t)=P_1U(t),~U_\infty(t)=P_\infty U(t).
\]
It follows that $U_1(t)=(\phi_1(t),w_1(t),\nabla\psi_1(t))$ and $U_\infty(t)=(\phi_\infty(t),w_\infty(t),\nabla\psi_\infty(t))$ satisfy the equations
\begin{gather}
\left\{
\begin{array}{l}
\displaystyle
U_1(t)=e^{-tL}U_1(0)+\int_{0}^{t} e^{-(t-s)L}P_1N(s)\mathrm{d}s, \\
\phi_1+\mathrm{div}\psi_1=0,
\end{array}
\right.
\label{lowfreq}  
\end{gather}
and
\begin{gather}
\left\{
\begin{array}{l}
\displaystyle
\partial_t U_\infty+LU_\infty=P_\infty N, \\
\phi_\infty+\mathrm{div}\psi_\infty=0.
\end{array}
\right.\label{highfreq}
\end{gather}

We first derive the $L^\infty$ estimate of the low frequency part $U_1(t)$.
\begin{prop}\label{LinftyestimateofUL}
There exists a positive number $\delta_0$ such that if $\|u_0\|_{L^1}+\|u_0\|_{H^3}\leq\delta_0$, then the following inequality 
$$
\|U_1(t)\|_{L^\infty}\leq C(1+t)^{-2}(\|u_0\|_{L^1}+\|u_0\|_{H^3})
$$
holds for $t\geq 0$.
\end{prop}
To prove Proposition \ref{LinftyestimateofUL}, we introduce the following lemmata.
\begin{lem}\label{L1Linftyestsemigroup}
Let $f\in L^1$. Then, the following estimates hold for $j\in\{0\}\cup\mathbb{N}$, $\alpha\in (\{0\}\cup\mathbb{N})^3$ and $t\geq0$
\begin{align*}
&\left\|\partial_t^j\partial_x^\alpha\mathcal{F}^{-1}\left[\frac{e^{\mu_1(\xi)t}-e^{\mu_2(\xi)t}}{\mu_1(\xi)-\mu_2(\xi)} \hat\eta(\xi)\hat\varphi_1(\xi)\right]\right\|_{L^\infty} \leq C(1+t)^{-\frac{3}{2}-\frac{j+|\alpha|}{2}}, 
\\ 
&\left\|\partial_t^j\partial_x^\alpha\mathcal{F}^{-1}\left[\frac{\mu_1(\xi)e^{\mu_2(\xi)t}-\mu_2(\xi)e^{\mu_1(\xi)t}}{\mu_1(\xi)-\mu_2(\xi)}\hat\eta(\xi)\hat\varphi_1(\xi)\right]\right\|_{L^\infty} \leq C(1+t)^{-2-\frac{j+|\alpha|}{2}}, 
\\ 
&\left\|\partial_t^j\partial_x^\alpha\mathcal{F}^{-1}\left[\frac{e^{\mu_3(\xi)t}-e^{\mu_4(\xi)t}}{\mu_3(\xi)-\mu_4(\xi)} \hat\eta(\xi)\hat\varphi_1(\xi)\right]\right\|_{L^\infty} \leq C(1+t)^{-\frac{3}{2}-\frac{j+|\alpha|}{2}}, 
\\ 
&\left\|\partial_t^j\partial_x^\alpha\mathcal{F}^{-1}\left[\frac{\mu_3(\xi)e^{\mu_4(\xi)t}-\mu_4(\xi)e^{\mu_3(\xi)t}}{\mu_3(\xi)-\mu_4(\xi)}\hat\eta(\xi)\hat\varphi_1(\xi)\right]\right\|_{L^\infty} \leq C(1+t)^{-2-\frac{j+|\alpha|}{2}}, 
\end{align*}
where $\hat\eta(\xi)=\tilde\eta(\frac{\xi}{|\xi|})$ with $\tilde\eta\in C^\infty(S^2)$ and $S^2=\{\omega\in\mathbb{R}^3||\omega|=1\}$. 
\end{lem}
Lemma \ref{L1Linftyestsemigroup} directly follows from \cite[Theorem3.1]{kobayashishibata}.

We give the estimate of $\|e^{-tL}U_1(0)\|_{L^\infty}$ as follows.
\begin{lem}\label{estetLu1}
It holds the following estimate:
\begin{gather*}
\|e^{-tL}U_1(0)\|_{L^\infty} \leq C(1+t)^{-2}\|u_0\|_{L^1}. 
\end{gather*}
\end{lem}
Lemma \ref{estetLu1} is a direct consequence of Lemma \ref{L1Linftyestsemigroup}.

For simplicity, we set $\|u_0\|_{\mathcal{X}}=\|u_0\|_{L^1}+\|u_0\|_{H^3}$. We have the estimate of $\int_{0}^{t} \|e^{-(t-s)L}P_1N(s)\|_{L^\infty}\mathrm{d}s$.
\begin{lem}\label{NLestLFLinf}
There exists a positive number $\delta_0$ such that if $\|u_0\|_{\mathcal{X}}\leq\delta_0$, then the following inequality 
\begin{align}\label{LinftyestUinfty}
\int_{0}^{t} \|e^{-(t-s)L}P_1N(s)\|_{L^\infty}\mathrm{d}s\leq C(1+t)^{-2}\|u_0\|_{\mathcal{X}}
\end{align}
holds for $t\geq 0$.
\end{lem}

\textbf{Proof.}
We first consider $\mathcal{F}^{-1}\left[\hat\varphi_1(\xi)\hat{K}^{23}(\xi,t-s)\hat{N}_3(\xi,s)\right]$ \\and $\mathcal{F}^{-1}\left[\hat\varphi_1(\xi)\hat{K}^{33}(\xi,t-s)\hat{N}_3(\xi,s)\right]$.
Since 
\[
\hat{N}_3(\xi,s)=-\left(I-\frac{\xi{}^\top\xi}{|\xi|^2}\right)\mathcal{F}(\nabla(w\cdot\nabla\tilde{\psi}))(\xi,s)+\frac{\xi{}^\top\xi}{|\xi|^2}\mathcal{F}(\nabla(\phi w))(\xi,s),
\]
we have
\begin{equation}
\begin{array}{l}
\displaystyle
\hat{K}^{23}(\xi,t-s)\hat{N}_3(\xi,s)\\[2ex]
\displaystyle
\quad=-i\beta^2\frac{e^{\mu_1(\xi)(t-s)}-e^{\mu_2(\xi)(t-s)}}{\mu_1(\xi)-\mu_2(\xi)}\left(I-\frac{{\xi}^\top \xi}{|\xi|^2}\right)\mathcal{F}(\nabla(w\cdot\nabla\tilde{\psi}))(\xi,s)\xi \\[2ex]
\displaystyle
\quad\quad\quad\quad+i\beta^2\frac{e^{\mu_3(\xi)(t-s)}-e^{\mu_4(\xi)(t-s)}}{\mu_3(\xi)-\mu_4(\xi)}\frac{{\xi}^\top \xi}{|\xi|^2}\mathcal{F}(\nabla(\phi w))(\xi,s)\xi, 
\end{array}
\label{calK23}
\end{equation}
\begin{equation}
\begin{array}{l}
\displaystyle
\hat{K}^{33}(\xi,t-s)\hat{N}_3(\xi,s)\\[2ex]
\displaystyle
\quad=\frac{\mu_1(\xi)e^{\mu_2(\xi)(t-s)}-\mu_2(\xi)e^{\mu_1(\xi)(t-s)}}{\mu_1(\xi)-\mu_2(\xi)}\left(I-\frac{{\xi}^\top \xi}{|\xi|^2}\right)\mathcal{F}(\nabla(w\cdot\nabla\tilde{\psi}))(\xi,s) \\[2ex]
\displaystyle
\quad\quad\quad\quad+\frac{\mu_3(\xi)e^{\mu_4(\xi)(t-s)}-\mu_4(\xi)e^{\mu_3(\xi)(t-s)}}{\mu_3(\xi)-\mu_4(\xi)}\frac{{\xi}^\top \xi}{|\xi|^2}\mathcal{F}(\nabla(\phi w))(\xi,s)), 
\end{array}
\label{calK33}
\end{equation}
We see from Lemma \ref{L1Linftyestsemigroup}, \eqref{calK23} and \eqref{calK33} that
\begin{equation}
\begin{array}{l}
\displaystyle
\left\|\mathcal{F}^{-1}\left[\hat\varphi_1(\xi)\hat{K}^{j3}(\xi,t-s)\hat{N}_3(\xi,s)\right]\right\|_{L^\infty}\\[2ex]
\displaystyle
\quad
\leq C(1+t-s)^{-2}(1+s)^{-2}\|u_0\|_{\mathcal{X}},~j=1,2,3. 
\end{array}
\label{estN33i} 
\end{equation}
Since
\begin{align}
\|N_1(s)\|_{L^1}&\leq C\|u(s)\|_{L^2}\|\nabla u(s)\|_{L^2}\leq C(1+s)^{-2}\|u_0\|_{\mathcal{X}}, \label{ineq5.12}\\
\|N_2(s)\|_{L^1}&\leq C\|u(s)\|_{L^2}\|\nabla u(s)\|_{H^1}\leq C(1+s)^{-2}\|u_0\|_{\mathcal{X}},\label{ineq5.13}
\end{align}
we see from Lemma \ref{L1Linftyestsemigroup} that 
\begin{equation}
\begin{array}{l}
\displaystyle
\left\|\mathcal{F}^{-1}\left[\hat\varphi_1(\xi)\hat{K}^{jk}(\xi,t-s)\hat{N}_k(\xi,s)\right]\right\|_{L^\infty}
\\[2ex]
\displaystyle
\quad \leq C(1+t-s)^{-2}(1+s)^{-2}\|u_0\|_{\mathcal{X}},~j=1,2,3,~k=1,2. 
\end{array}
\label{estN32i} 
\end{equation}
It follows from \eqref{estN33i} and \eqref{estN32i} that
\[
\|e^{-(t-s)L}P_1N(s)\|_{L^\infty}\leq C(1+t-s)^{-2}(1+s)^{-2}\|u_0\|_{\mathcal{X}}.
\] 
By employing Lemma \ref{betafunc1} with $a=b=2$, we have \eqref{LinftyestUinfty}. 
This completes the proof. $\blacksquare$

\textbf{Proof of Proposition \ref{LinftyestimateofUL}.}
Taking $L^\infty$ norm of the first equation of \eqref{lowfreq}, we have
\begin{align}\label{lowfreqineq}
\|U_1(t)\|_{L^\infty}&\leq \|e^{-tL}U_1(0)\|_{L^\infty}+\int_{0}^{t} \|e^{-(t-s)L}P_1N(s)\|_{L^\infty}\mathrm{d}s. 
\end{align}
Together with \eqref{lowfreqineq}, Lemma \ref{estetLu1} and Lemma \ref{NLestLFLinf}, then yields
\begin{align}
\|U_1(t)\|_{L^\infty}\leq C(1+t)^{-2}\|u_0\|_{\mathcal{X}}.
\label{LinftyestU1} 
\end{align}
This completes the proof of Proposition \ref{LinftyestimateofUL}. $\blacksquare$

We next consider the high frequency part $U_\infty(t)$.
\begin{prop}\label{LinftyestimateofUinfty}
There exists a positive number $\delta_0$ such that if $\|u_0\|_{L^1}+\|u_0\|_{H^3}\leq\delta_0$, then the following inequality
$$
\|U_\infty(t)\|_{L^\infty}\leq C(1+t)^{-2}(\|u_0\|_{L^1}+\|u_0\|_{H^3})
$$
holds for $t\geq 0$.
\end{prop}
\textbf{Proof. }
We set $\tilde{\psi}_\infty=P_\infty\tilde{\psi}$, $\tilde{\varPsi}_\infty=\nabla\tilde{\psi}_\infty$ and $\tilde{U}_\infty=(\phi_\infty,w_\infty,\nabla\tilde{\psi}_\infty)$.

Since 
\begin{align*}
\|\varPsi_\infty\|_{L^\infty}
&\leq c\|\varPsi_\infty\|_{H^2} \\
&\leq c\|\tilde{\varPsi}_\infty\|_{H^2}+c\| P_\infty\nabla(-\Delta)^{-1}\mathrm{div}{}^\top(\phi\nabla\tilde{\psi}+(1+\phi)h(\nabla\tilde{\psi}))\|_{H^2} \\
&\leq c\|\tilde{\varPsi}_\infty\|_{H^2}+C\|\nabla P_\infty\nabla(-\Delta)^{-1}\mathrm{div}{}^\top(\phi\nabla\tilde{\psi}+(1+\phi)h(\nabla\tilde{\psi}))\|_{H^1} \\
&\leq c\|\tilde{\varPsi}_\infty\|_{H^2}+C(1+t)^{-2}\|u_0\|_{\mathcal{X}}, 
\end{align*}
 we have $\|U_\infty(t)\|_{L^\infty}\leq c\|\tilde{U}_\infty(t)\|_{H^2}+C(1+t)^{-2}\|u_0\|_{\mathcal{X}}$. We thus estimate $\tilde{\psi}_\infty$ in substitution for $\psi_\infty$.

We next consider $\|\tilde{U}_\infty\|_{L^\infty}$ instead of $\|U_\infty\|_{L^\infty}$. 
By applying $P_\infty$ to the problem \eqref{system3}, we obtain
\begin{gather}
\partial_t\phi_\infty+\mathrm{div} w_\infty=f_{1,\infty}, \label{HFeq1}\\
\partial_t  w_\infty-\nu{\Delta} w_\infty-\tilde \nu\nabla \mathrm{div}  w_\infty+\gamma^2\nabla\phi_\infty-\beta^2\mathrm{div} \tilde{\varPsi}_\infty = f_{2,\infty}, \label{HFeq2} \\
\partial_t \tilde{\varPsi}_\infty-\nabla w_\infty= f_{3,\infty}, \label{HFeq3} \\
\phi_\infty  +\mathrm{div}\tilde{\psi}_\infty=\tilde{f}_{4,\infty}, \label{HFeq4} \\
\tilde{U}_\infty|_{t=0}=P_\infty(\phi_0,w_0,\tilde{\varPsi}_0)\nonumber,
\end{gather}
where $f_{j,\infty}=P_\infty f_j,~j=1,2,3$, and $\tilde{f}_{4,\infty}=-P_\infty (-\Delta)^{-1}\mathrm{div}f_4$.
We define $E[\tilde{U}_\infty]$ and $D[\tilde{U}_\infty]$ by
$$
E[\tilde{U}_\infty]=\|\tilde{U}_\infty\|^{2}_{H^2}+c_1\sum_{|\alpha|\leq2}(\partial_x^\alpha w_\infty,\partial_x^\alpha\tilde{\psi}_\infty),
$$ 
$$
D[\tilde{U}_\infty]=\sum_{|\alpha|\leq2}[\nu\|\nabla\partial_x^\alpha w_\infty\|^{2}_{L^2}+\tilde{\nu}\|\mathrm{div}\partial_x^\alpha w_\infty\|^{2}_{L^2}+c_1\gamma^2\|\partial_x^\alpha \phi_\infty\|^{2}_{L^2}+c_1\beta^2\|\partial_x^\alpha \tilde{\varPsi}_\infty\|^{2}_{L^2}].
$$
Here $c_1$ is a positive constant to be determined later.

We establish the following energy estimate of $E[\tilde{U}_\infty]$.

\begin{prop}\label{estforHF}
The following estimate holds:
\begin{equation}
\frac{\mathrm{d}}{\mathrm{d}t}E[\tilde{U}_\infty]+D[\tilde{U}_\infty]\leq C_1\mathcal{N}. \label{energyestforUinfty}
\end{equation}
Here
\begin{align*} 
\mathcal{N}&=\sum_{|\alpha|\leq2}\left(\gamma^2|(\partial_x^\alpha f_{1,\infty},\partial_x^\alpha \phi_\infty)|+|(\partial_x^\alpha f_{2,\infty},\partial_x^\alpha w_\infty)|+\beta^2|(\partial_x^\alpha f_{3,\infty},\partial_x^\alpha\tilde{\varPsi}_\infty)|\right. \\
&\quad\quad+c_1|(\partial_x^\alpha f_{2,\infty},\partial_x^\alpha \tilde{\psi}_\infty)|+c_1|(\partial_x^\alpha (-\Delta)^{-1}\mathrm{div}f_{3,\infty},\partial_x^\alpha w_\infty)| \\[1ex]
&\quad\quad\quad\left.+c_1\gamma^2|(\nabla\partial_x^\alpha \tilde{f}_{4,\infty},\partial_x^\alpha \tilde{\psi}_\infty)|+c_1\gamma^2\|\partial_x^\alpha \tilde{f}_{4,\infty}\|_{L^2}^2
\right).
\end{align*}
\end{prop}
\textbf{Proof.} We take the inner product of $\partial_x^\alpha\eqref{HFeq1}$ with $\partial_x^\alpha\phi_\infty$ to obtain
\begin{gather}
\frac{1}{2}\frac{\mathrm{d}}{\mathrm{d}t}\|\partial_x^\alpha\phi_\infty\|^{2}_{L^2}-(\nabla\partial_x^\alpha\phi_\infty,\partial_x^\alpha w_\infty)=(\partial_x^\alpha f_{1,\infty}, \partial_x^\alpha\phi_\infty). \label{eeq1}
\end{gather}
We take the inner product of $\partial_x^\alpha\eqref{HFeq2}$ with $\partial_x^\alpha w_\infty$ to obtain
\begin{gather}
\begin{array}{l}
\displaystyle
\frac{1}{2}\frac{\mathrm{d}}{\mathrm{d}t}\|\partial_x^\alpha w_\infty\|^{2}_{L^2}+\nu\|\nabla \partial_x^\alpha w_\infty\|^{2}_{L^2}
+\tilde{\nu}\|\mathrm{div}\partial_x^\alpha w_\infty\|^{2}_{L^2}\\[2ex]
\displaystyle
+\gamma^2(\nabla\partial_x^\alpha \phi_\infty,\partial_x^\alpha w_\infty)
-\beta^2(\mathrm{div} \partial_x^\alpha\tilde{\varPsi}_\infty,\partial_x^\alpha w_\infty)
= (\partial_x^\alpha f_{2,\infty},\partial_x^\alpha w_\infty).
\end{array}\label{eeq2}
\end{gather}
We take the inner product of $\partial_x^\alpha\eqref{HFeq3}$ with $\partial_x^\alpha\tilde{\varPsi}_\infty$ to obtain
\begin{gather}
\frac{1}{2}\frac{\mathrm{d}}{\mathrm{d}t}\|\partial_x^\alpha\tilde{\varPsi}_\infty\|^{2}_{L^2}+(\mathrm{div} \partial_x^\alpha\tilde{\varPsi}_\infty,\partial_x^\alpha w_\infty) = (\partial_x^\alpha f_{3,\infty},\partial_x^\alpha\tilde{\varPsi}_\infty).\label{eeq3}
\end{gather}
It then follows from $\gamma^2\times(\ref{eeq1})+(\ref{eeq2})+\beta^2\times(\ref{eeq3})$ that
\begin{gather}
\begin{array}{l}
\displaystyle{\frac{1}{2}\frac{\mathrm{d}}{\mathrm{d}t}}(\gamma^2\|\partial_x^\alpha\phi_\infty\|^{2}_{L^2}+\|\partial_x^\alpha w_\infty\|^{2}_{L^2}+\beta^2\|\partial_x^\alpha\tilde{\varPsi}_\infty\|^{2}_{L^2})+D^0[\partial_x^\alpha w_\infty] \\[2ex]
\displaystyle
\quad
\leq \gamma^2|(\partial_x^\alpha f_{1,\infty},\partial_x^\alpha \phi_\infty)|+|(\partial_x^\alpha f_{2,\infty},\partial_x^\alpha w_\infty)|+\beta^2|(\partial_x^\alpha f_{3,\infty},\partial_x^\alpha\tilde{\varPsi}_\infty)|.
\end{array}
\label{EIHF1}
\end{gather} 
Here $D^0[\partial_x^\alpha w_\infty]$ is given by
$$
D^0[\partial_x^\alpha w_\infty]=\nu\|\nabla\partial_x^\alpha w_\infty\|^{2}_{L^2}+\Tilde{\nu}\|\mathrm{div}\partial_x^\alpha w_\infty\|^{2}_{L^2}.
$$

We next derive the dissipative estimate of $\|\partial_x^\alpha\tilde{\varPsi}_\infty\|_{L^2}$. 
By substituting \eqref{HFeq4} to \eqref{HFeq2}, we have
\begin{equation}
\partial_t  w_\infty-\nu{\Delta} w_\infty-\tilde \nu\nabla \mathrm{div}  w_\infty-\beta^2\Delta\tilde{\psi}_\infty-\gamma^2\nabla\mathrm{div} \tilde{\psi}_\infty= f_{2,\infty}-\gamma^2\nabla\tilde{f}_{4,\infty}. \label{HFeq5}
\end{equation}
We take the inner product of $\partial_x^\alpha\eqref{HFeq5}$ with $\partial_x^\alpha\tilde{\psi}_\infty$ to obtain
\begin{gather}
\begin{array}{l}
(\partial_t\partial_x^\alpha w_\infty,\partial_x^\alpha\tilde{\psi}_\infty)-\nu({\Delta}\partial_x^\alpha w_\infty,\partial_x^\alpha\tilde{\psi}_\infty)
-\tilde{\nu}(\nabla \mathrm{div}\partial_x^\alpha  w_\infty,\partial_x^\alpha\tilde{\psi}_\infty) \\[1ex]
\quad
-\beta^2({\Delta}\partial_x^\alpha\tilde{\psi}_\infty,\partial_x^\alpha\tilde{\psi}_\infty)-\gamma^2(\nabla\mathrm{div}\partial_x^\alpha\tilde{\psi}_\infty,\partial_x^\alpha\tilde{\psi}_\infty)\\[1ex]
\quad\quad
=(\partial_x^\alpha f_{2,\infty},\partial_x^\alpha\tilde{\psi}_\infty)-\gamma^2(\nabla\partial_x^\alpha \tilde{f}_{4,\infty},\partial_x^\alpha\tilde{\psi}_\infty). 
\end{array}
\label{eq6.3.1}
\end{gather}
The first term on the left-hand side of $(\ref{eq6.3.1})$ is written as
\begin{align*}
&(\partial_t\partial_x^\alpha w,\partial_x^\alpha\tilde{\psi}_\infty)\\
&=-\frac{\mathrm{d}}{\mathrm{d}t}(\partial_x^\alpha w,\partial_x^\alpha\tilde{\psi}_\infty)+(\partial_x^\alpha w_\infty,\partial_t\partial_x^\alpha\tilde{\psi}_\infty) \\
&=-\frac{\mathrm{d}}{\mathrm{d}t}(\partial_x^\alpha w_\infty,\partial_x^\alpha\tilde{\psi}_\infty)-\|\partial_x^\alpha w_\infty\|_{L^2}^2-(\partial_x^\alpha (-\Delta)^{-1}\mathrm{div}f_{3,\infty},\partial_x^\alpha w_\infty).
\end{align*}
By integration by parts, the fourth term and fifth term of $(\ref{eq6.3.1})$ are written as
$-({\Delta}\partial_x^\alpha\tilde{\psi}_\infty,\partial_x^\alpha\tilde{\psi}_\infty)=\|\nabla \partial_x^\alpha\tilde{\psi}_\infty\|^{2}_{L^2}$ and $-(\nabla\mathrm{div}\partial_x^\alpha \tilde{\psi}_\infty,\partial_x^\alpha\tilde{\psi}_\infty)=\|\mathrm{div}\partial_x^\alpha \tilde{\psi}_\infty\|_{L^2}^2$, respectively.
It then follows from $(\ref{eq6.3.1})$ that
\begin{gather}
\begin{array}{l}
\displaystyle
-\frac{\mathrm{d}}{\mathrm{d}t}(\partial_x^\alpha w_\infty,\partial_x^\alpha\tilde{\psi}_\infty)+\beta^2\|\nabla \partial_x^\alpha\tilde{\psi}_\infty\|^{2}_{L^2}
+\gamma^2\|\mathrm{div}\partial_x^\alpha \tilde{\psi}_\infty\|_{L^2}^2\\[2ex]
\displaystyle
=\nu({\Delta}\partial_x^\alpha w_\infty,\partial_x^\alpha\tilde{\psi}_\infty)
+\tilde{\nu}(\nabla \mathrm{div}\partial_x^\alpha  w_\infty,\partial_x^\alpha\tilde{\psi}_\infty)
+\|\partial_x^\alpha w_\infty\|_{L^2}^2\\[2ex]
\displaystyle
+(\partial_x^\alpha f_{2,\infty},\partial_x^\alpha\tilde{\psi}_\infty)+(\partial_x^\alpha (-\Delta)^{-1}\mathrm{div}f_{3,\infty},\partial_x^\alpha w_\infty)-\gamma^2(\nabla\partial_x^\alpha \tilde{f}_{4,\infty},\partial_x^\alpha\tilde{\psi}_\infty). 
\end{array}
\label{eeq4}
\end{gather}
By integration by parts and the inequality $\|\partial_x^\alpha w_\infty\|_{L^2}\leq\|\nabla\partial_x^\alpha w_\infty\|_{L^2}$ following from Lemma \ref{lemP1Pinfty} (iii), we have 
\begin{align*}
 &\nu({\Delta}\partial_x^\alpha w_\infty,\partial_x^\alpha\tilde{\psi}_\infty)
+\tilde{\nu}(\nabla \mathrm{div}\partial_x^\alpha  w_\infty,\partial_x^\alpha\tilde{\psi}_\infty)
+\|\partial_x^\alpha w_\infty\|_{L^2}^2\\
&=-\nu(\nabla\partial_x^\alpha w_\infty,\nabla\partial_x^\alpha\tilde{\psi}_\infty)
+\tilde{\nu}( \mathrm{div}\partial_x^\alpha  w_\infty,\mathrm{div}\partial_x^\alpha\tilde{\psi}_\infty)
+\|\partial_x^\alpha w_\infty\|_{L^2}^2\\
&\leq \left(\frac{\nu}{2\beta^2}+\frac{\tilde{\nu}}{2\gamma^2}+\frac{1}{\nu}\right)D^0[\partial_x^\alpha w_\infty]+\frac{\beta^2}{2}\|\nabla \partial_x^\alpha\tilde{\psi}_\infty\|^{2}_{L^2}
+\frac{\gamma^2}{2}\|\mathrm{div}\partial_x^\alpha \tilde{\psi}_\infty\|_{L^2}^2.
\end{align*}
It then follows from $(\ref{eeq4})$ that
\begin{gather}
\begin{array}{l}
\displaystyle
-2\frac{\mathrm{d}}{\mathrm{d}t}(\partial_x^\alpha w_\infty,\partial_x^\alpha\tilde{\psi}_\infty)+\beta^2\|\nabla \partial_x^\alpha\tilde{\psi}_\infty\|^{2}_{L^2}
+\gamma^2\|\mathrm{div}\partial_x^\alpha \tilde{\psi}_\infty\|_{L^2}^2\\[2ex]
\displaystyle
\quad\leq\left(\frac{\nu}{\beta^2}+\frac{\tilde{\nu}}{\gamma^2}+\frac{2}{\nu}\right)D^0[\partial_x^\alpha w_\infty]+|(\partial_x^\alpha f_{2,\infty},\partial_x^\alpha\tilde{\psi}_\infty)|\\[2ex]
\displaystyle
\quad\quad+|(\partial_x^\alpha (-\Delta)^{-1}\mathrm{div}f_{3,\infty},\partial_x^\alpha w_\infty)|+\gamma^2|(\nabla\partial_x^\alpha \tilde{f}_{4,\infty},\partial_x^\alpha\tilde{\psi}_\infty)|. 
\end{array}
\label{EIHF2}
\end{gather}
Adding \eqref{EIHF1} to $\frac{c_1}{2}\times\eqref{EIHF2}$ and using $\|\partial_x^\alpha \phi_\infty\|_{L^2}^2\leq 2(\|\mathrm{div}\partial_x^\alpha \tilde{\psi}_\infty\|_{L^2}^2+\|\partial_x^\alpha \tilde{f}_{4,\infty}\|_{L^2}^2)$, we have
\begin{gather}
\begin{array}{l}
\displaystyle{\frac{1}{2}\frac{\mathrm{d}}{\mathrm{d}t}}(\gamma^2\|\partial_x^\alpha\phi_\infty\|^{2}_{L^2}+\|\partial_x^\alpha w_\infty\|^{2}_{L^2}+\beta^2\|\partial_x^\alpha\tilde{\varPsi}_\infty\|^{2}_{L^2}-2c_1(\partial_x^\alpha w_\infty,\partial_x^\alpha\tilde{\psi}_\infty)) \\[2ex]
\displaystyle
+\left(1-c_1\left(\frac{\nu}{\beta^2}+\frac{\tilde{\nu}}{\gamma^2}+\frac{2}{\nu}\right)\right)D^0[\partial_x^\alpha w_\infty]+\frac{c_1}{2}\gamma^2\|\partial_x^\alpha \phi_\infty\|_{L^2}^2+c_1\beta^2\| \partial_x^\alpha\tilde{\varPsi}_\infty\|^{2}_{L^2}
\\[2ex]
\displaystyle
\quad
\leq \gamma^2|(\partial_x^\alpha f_{1,\infty},\partial_x^\alpha \phi_\infty)|+|(\partial_x^\alpha f_{2,\infty},\partial_x^\alpha w_\infty)|+\beta^2|(\partial_x^\alpha f_{3,\infty},\partial_x^\alpha\tilde{\varPsi}_\infty)| \\[2ex]
\displaystyle
\quad\quad+2c_1|(\partial_x^\alpha f_{2,\infty},\partial_x^\alpha\tilde{\psi}_\infty)|+2c_1|(\partial_x^\alpha (-\Delta)^{-1}\mathrm{div}f_{3,\infty},\partial_x^\alpha w_\infty)|\\[2ex]
\displaystyle
\quad\quad\quad+2c_1\gamma^2|(\nabla\partial_x^\alpha {f}_{4,\infty},\partial_x^\alpha\tilde{\psi}_\infty)|+c_1\gamma^2\|\partial_x^\alpha \tilde{f}_{4,\infty}\|_{L^2}^2. 
\end{array}
\label{EEHF3}
\end{gather}
We take $c_1>0$ small so that $1-c_1\left(\frac{\nu}{\beta^2}+\frac{\tilde{\nu}}{\gamma^2}+\frac{2}{\nu}\right)\geq\frac{1}{2}$ and $d_1\|\tilde{U}_\infty(t)\|_{H^2}^2\leq E[\tilde{U}_\infty](t)\leq d_2D[\tilde{U}_\infty](t)$ for some positive numbers $d_1,~d_2>0$. By summing \eqref{EEHF3} for $|\alpha|\leq2$, we obtain \eqref{energyestforUinfty}.
This completes the proof. $\blacksquare$

We next estimate $\mathcal{N}(t)$. 
\begin{prop}\label{estofmathcalN}
The following estimate holds uniformly in $t\geq0$
\begin{align}
\mathcal{N}(t)\leq C_2\delta_0 D[\tilde{U}_\infty](t)+C_2(1+t)^{-4}\|u_0\|_{\mathcal{X}}^2. \label{NLestHF1i}
\end{align}
\end{prop}
Proposition \ref{estofmathcalN} can be shown by using Lemma \ref{Sobolevineq}, Proposition \ref{H2estimate}, and integration by parts. We omit the proof.

\textbf{Proof of Proposition \ref{LinftyestimateofUinfty} (continued). }
By taking $\delta_0$ so that $C_1C_2\delta_0\leq\frac{1}{2}$, it follows from \eqref{energyestforUinfty} and \eqref{NLestHF1i} that
\begin{align}
\frac{\mathrm{d}}{\mathrm{d}t}E[\tilde{U}_\infty]+\frac{1}{2d_2}E[\tilde{U}_\infty]\leq C(1+t)^{-4}\|u_0\|_{\mathcal{X}}^2. \label{energyestforUinfty2}
\end{align}
Therefore, we see from Lemma \ref{betafunc1} with $a=b=4$ and \eqref{energyestforUinfty2} that 
\begin{align*}
E[\tilde{U}_\infty](t) &\leq e^{-ct}E[\tilde{U}_\infty](0)+C\int_0^t e^{-c(t-s)}(1+s)^{-4}\mathrm{d}s \|u_0\|_{\mathcal{X}}^2\\
             		 &\leq C(1+t)^{-4}\|\tilde{U}_0\|_{H^2}^2+C\int_0^t (1+t-s)^{-4}(1+s)^{-4}\mathrm{d}s\|u_0\|_{\mathcal{X}}^2\\
                   &\leq C(1+t)^{-4}\|u_0\|_{\mathcal{X}}^2. 
\end{align*}
Since $\|U_\infty(t)\|_{L^\infty}^2\leq CE[\tilde{U}_\infty](t)+C(1+t)^{-4}\|u_0\|_{\mathcal{X}}^2$, we finally arrive at
\begin{align*}
\|U_\infty(t)\|_{L^\infty}\leq C(1+t)^{-2}\|u_0\|_{\mathcal{X}}. 
\end{align*}
This completes the proof of Proposition \ref{LinftyestimateofUinfty}. $\blacksquare$

\textbf{Proof of Proposition \ref{LinftyestimateofU}.}  Proposition \ref{LinftyestimateofU} immediately follows from Proposition \ref{LinftyestimateofUL} and Proposition \ref{LinftyestimateofUinfty}. This completes the proof. $\blacksquare$

\section{Proof of Theorem \ref{mainthm} (ii)} 
In this section, we give a proof of Theorem \ref{mainthm} (ii).
In view of Lemmata \ref{equivforGpsi}--\ref{psitildepsi},
it suffices to obtain the following $L^p$ decay estimate of $U=(\phi,w,\varPsi)$ to prove Theorem \ref{mainthm} (ii).
\begin{prop}\label{L1estimate}
There exists a positive number $\delta_0$ such that if $\|u_0\|_{L^1}+\|u_0\|_{H^3}\leq\delta_0$, then the following inequality 
$$
\|U(t)\|_{L^p}\leq C(p)(1+t)^{-\frac{3}{2}\left(1-\frac{1}{p}\right)+\frac{1}{2}\left(\frac{2}{p}-1\right)}(\|u_0\|_{L^1}+\|u_0\|_{L^p}+\|u_0\|_{H^3})
$$
holds for $1<p<2$ and $t\geq 0$.
\end{prop}
Proposition \ref{L1estimate} is a direct consequence of the $L^p$ estimates of $U_1(t)$ and $U_\infty(t)$ which will be established below in Proposition \ref{L1Lplowfreq} and Proposition \ref{LpLphighfreq}, respectively.

We first consider the low frequency part $U_1(t)$.

\begin{prop}\label{L1Lplowfreq}
There exists a positive number $\delta_0$ such that if $\|u_0\|_{\mathcal{X}}\leq\delta_0$, then it holds the following estimate:
\begin{align*}
\|U_1(t)\|_{L^p}\leq C(1+t)^{-\frac{3}{2}\left(1-\frac{1}{p}\right)+\frac{1}{2}\left(\frac{2}{p}-1\right)}\|u_0\|_{\mathcal{X}}.
\end{align*}
\end{prop}
Since $\|U_1(t)\|_{L^p}\leq\|U_1(t)\|_{L^1}^{\frac{2}{p}-1}\|U_1(t)\|_{L^2}^{2-\frac{2}{p}},~1< p<2,$ it is enough to show the case $p=1$ only:
\begin{align}
\|U_1(t)\|_{L^1}\leq C(1+t)^{\frac{1}{2}}\|u_0\|_{\mathcal{X}}. \label{L1L1lowfreq}
\end{align}
To show \eqref{L1L1lowfreq}, we introduce the following lemma.
\begin{lem}\label{L1L1LFKS}
Let $f\in L^1$. Then, the following estimates hold for $j\geq0$, $\alpha\in (\{0\}\cup\mathbb{N})^3$ and $t\geq0$:
\begin{align*}
&\left\|\partial_t^j\partial_x^\alpha\mathcal{F}^{-1}\left[\frac{e^{\mu_1(\xi)t}-e^{\mu_2(\xi)t}}{\mu_1(\xi)-\mu_2(\xi)} \hat\varphi_1(\xi)\right]\right\|_{L^1} \leq C(1+t)^{1-\frac{j+|\alpha|}{2}}, 
\\ 
&\left\|\partial_t^j\partial_x^\alpha\mathcal{F}^{-1}\left[\frac{\mu_1(\xi)e^{\mu_2(\xi)t}-\mu_2(\xi)e^{\mu_1(\xi)t}}{\mu_1(\xi)-\mu_2(\xi)}\hat\varphi_1(\xi)\right]\right\|_{L^1} \leq C(1+t)^{\frac{1}{2}-\frac{j+|\alpha|}{2}}, 
\\ 
&\left\|\partial_t^j\partial_x^\alpha\mathcal{F}^{-1}\left[\frac{e^{\mu_3(\xi)t}-e^{\mu_4(\xi)t}}{\mu_3(\xi)-\mu_4(\xi)} \hat\varphi_1(\xi)\right]\right\|_{L^1} \leq C(1+t)^{1-\frac{j+|\alpha|}{2}}, 
\\ 
&\left\|\partial_t^j\partial_x^\alpha\mathcal{F}^{-1}\left[\frac{\mu_3(\xi)e^{\mu_4(\xi)t}-\mu_4(\xi)e^{\mu_3(\xi)t}}{\mu_3(\xi)-\mu_4(\xi)}\hat\varphi_1(\xi)\right]\right\|_{L^1} \leq C(1+t)^{\frac{1}{2}-\frac{j+|\alpha|}{2}}, 
\\
&\left\|\partial_t^j\partial_x^\alpha\mathcal{F}^{-1}\left[\left(\frac{\mu_1(\xi)e^{\mu_1(\xi)t}-\mu_2(\xi)e^{\mu_2(\xi)t}}{\mu_1(\xi)-\mu_2(\xi)}-e^{-\nu|\xi|^2t}\right)\frac{\xi{}^\top\xi}{|\xi|}\hat\varphi_1(\xi)\right]\right\|_{L^1}\nonumber \\
&\leq C(1+t)^{\frac{1}{2}-\frac{j+|\alpha|}{2}}, 
\\
&\left\|\partial_t^j\partial_x^\alpha\mathcal{F}^{-1}\left[\left(\frac{\mu_3(\xi)e^{\mu_3(\xi)t}-\mu_4(\xi)e^{\mu_4(\xi)t}}{\mu_3(\xi)-\mu_4(\xi)}-e^{-\nu|\xi|^2t}\right)\frac{\xi{}^\top\xi}{|\xi|}\hat\varphi_1(\xi)\right]\right\|_{L^1} 
\nonumber \\
&\leq C(1+t)^{\frac{1}{2}-\frac{j+|\alpha|}{2}}. 
\end{align*}
\end{lem}
Lemma \ref{L1L1LFKS} is obtained in \cite[pp.216]{kobayashishibata} and \cite[pp.216]{shibataVW2000} directly.

We have the estimate of $\|e^{-tL}U_1(0)\|_{L^1}$.
\begin{lem}\label{L1L1est-low-semigroup}
The following estimate holds for $t\geq0$:
\[
\|e^{-tL} U_1(0)\|_{L^1}\leq C(1+t)^{\frac{1}{2}}\|u_0\|_{L^1}.
\]
\end{lem}

\textbf{Proof.}
The $L^1$ estimates of $\mathcal{F}^{-1}[\hat\varphi_1(\xi)\hat{K}^{j1}(\xi,t)\hat\phi_0(\xi)]$ $(j=1,2,3)$, and
$\mathcal{F}^{-1}[\hat\varphi_1(\xi)\hat{K}^{12}(\xi,t)\hat{w}_0(\xi)]$ immediately follow from Lemma \ref{L1L1LFKS}:
\begin{align}
\|\mathcal{F}^{-1}[\hat\varphi_1(\xi)\hat{K}^{j1}(\xi,t)\hat\phi_0(\xi)]\|_{L^1}&\leq C(1+t)^{\frac{1}{2}}\|\phi_0\|_{L^1},~j=1,2,3, \label{ineq6.7}\\
\|\mathcal{F}^{-1}[\hat\varphi_1(\xi)\hat{K}^{12}(\xi,t)\hat{w}_0(\xi)]\|_{L^1}  &\leq C(1+t)^{\frac{1}{2}}\|w_0\|_{L^1}. \label{ineq6.9} 
\end{align}
Since
\begin{align*}
\hat{K}^{22}(\xi,t)\hat{w}_0&=\frac{\mu_1(\xi)e^{\mu_1(\xi)t}-\mu_2(\xi)e^{\mu_2(\xi)t}}{\mu_1(\xi)-\mu_2(\xi)}\hat{w}_0(\xi) \\
&-\left(\frac{\mu_1(\xi)e^{\mu_1(\xi)t}-\mu_2(\xi)e^{\mu_2(\xi)t}}{\mu_1(\xi)-\mu_2(\xi)}-e^{-\nu|\xi|^2t}\right)\frac{{\xi}^\top \xi}{|\xi|^2}\hat{w}_0(\xi)\\
&\quad+\left(\frac{\mu_3(\xi)e^{\mu_3(\xi)t}-\mu_4(\xi)e^{\mu_4(\xi)t}}{\mu_3(\xi)-\mu_4(\xi)}-e^{-\nu|\xi|^2t}\right)\frac{{\xi}^\top \xi}{|\xi|^2}\hat{w}_0(\xi),
\end{align*}
\begin{align*}
\hat{K}^{33}(\xi,t)\hat\varPsi_0&=\frac{\mu_1(\xi)e^{\mu_2(\xi)t}-\mu_2(\xi)e^{\mu_1(\xi)t}}{\mu_1(\xi)-\mu_2(\xi)}\hat{\varPsi}_0(\xi) \\
&-\left(\frac{\mu_1(\xi)e^{\mu_1(\xi)t}-\mu_2(\xi)e^{\mu_2(\xi)t}}{\mu_1(\xi)-\mu_2(\xi)}-e^{-\nu|\xi|^2t}\right)\frac{{\xi}^\top \xi}{|\xi|^2}\hat{\varPsi}_0(\xi) \\
&\quad-\nu\frac{e^{\mu_1(\xi)t}-e^{\mu_2(\xi)t}}{\mu_1(\xi)-\mu_2(\xi)}\xi^\top \xi\hat\varPsi_0(\xi)\\ 
&+\left(\frac{\mu_3(\xi)e^{\mu_3(\xi)t}-\mu_4(\xi)e^{\mu_4(\xi)t}}{\mu_3(\xi)-\mu_4(\xi)}-e^{-\nu|\xi|^2t}\right)\frac{{\xi}^\top \xi}{|\xi|^2}\hat{\varPsi}_0(\xi)\\
&\quad+(\nu+\tilde{\nu})\frac{e^{\mu_3(\xi)t}-e^{\mu_4(\xi)t}}{\mu_3(\xi)-\mu_4(\xi)}\xi^\top \xi\hat\varPsi_0(\xi) ,
\end{align*}
we see from Lemma \ref{L1L1LFKS} that
\begin{align}
\|\mathcal{F}^{-1}[\hat\varphi_1(\xi)\hat{K}^{22}(\xi,t)\hat{w}_0(\xi)]\|_{L^1}&\leq C(1+t)^{\frac{1}{2}}\|w_0\|_{L^1}, \label{ineq6.10}\\
\|\mathcal{F}^{-1}[\hat\varphi_1(\xi)\hat{K}^{33}(\xi,t)\hat\varPsi_0(\xi)]\|_{L^1}&\leq C(1+t)^{\frac{1}{2}}\|\varPsi_0\|_{L^1}. \label{ineq6.11} 
\end{align}
It remains to estimate $\mathcal{F}^{-1}[\hat\varphi_1(\xi)\hat{K}^{23}(\xi,t)\hat\varPsi_0]$ and $\mathcal{F}^{-1}[\hat\varphi_1(\xi)\hat{K}^{32}(\xi,t)\hat{w}_0]$.

We write $(\mathcal{F}^{-1}[\hat\varphi_1(\xi)\hat{K}^{23}(\xi,t)\hat\varPsi_0])^j,~j=1,2,3,$ and $(\mathcal{F}^{-1}[\hat\varphi_1(\xi)\hat{K}^{32}(\xi,t)\hat{w}_0])^{jk},$ $j,k=1,2,3,$ as
\begin{align*}
&(\mathcal{F}^{-1}[\hat\varphi_1(\xi)\hat{K}^{23}(\xi,t)\hat\varPsi_0])^j \\
&=\beta^2\sum_{k=1}^3(\mathcal{K}_0^{k}\ast\hat{\varPsi}_0^{jk})(\xi) 
-\beta^2\sum_{k,l=1}^3(\mathcal{K}_1\ast\mathcal{L}_{j,k,l}^{\frac{\nu}{4}}\ast\hat{\varPsi}_0^{lk})(\xi) \\
&\quad+\beta^2\sum_{k,l=1}^3(\mathcal{K}_2\ast\mathcal{L}_{j,k,l}^{\frac{\nu+\tilde{\nu}}{4}}\ast\hat{\varPsi}_0^{lk})(\xi) ,\\
&(\mathcal{F}^{-1}[\hat\varphi_1(\xi)\hat{K}^{32}(\xi,t)\hat{w}_0])^{jk}\\
&=(\mathcal{K}_0^{k}\ast\hat{w}_0^{j})(\xi)  
-\sum_{l=1}^3(\mathcal{K}_1\ast\mathcal{L}_{j,k,l}^{\frac{\nu}{4}}\ast\hat{w}_0^{l})(\xi)
+\sum_{l=1}^3(\mathcal{K}_2\ast\mathcal{L}_{j,k,l}^{\frac{\nu+\tilde{\nu}}{4}}\ast\hat{w}_0^{l})(\xi) ,
\end{align*}
where
\begin{align*}
\mathcal{K}_0^{k}&=\mathcal{F}^{-1}\left[i\xi_k\frac{e^{\mu_1(\xi)t}-e^{\mu_2(\xi)t}}{\mu_1(\xi)-\mu_2(\xi)}\hat\varphi_1(\xi)\right],~k=1,2,3, \\
\mathcal{K}_1&=\mathcal{F}^{-1}\left[\frac{e^{\tilde\mu_1(\xi)t}-e^{\tilde\mu_2(\xi)t}}{\tilde\mu_1(\xi)-\tilde\mu_2(\xi)}\hat\varphi_1(\xi)\right],\\
\mathcal{K}_2&=\mathcal{F}^{-1}\left[\frac{e^{\tilde\mu_3(\xi)t}-e^{\tilde\mu_4(\xi)t}}{\tilde\mu_3(\xi)-\tilde\mu_4(\xi)}\hat\varphi_1(\xi)\right],\\
\mathcal{L}_{l,jk}^a&=\mathcal{F}^{-1}\left[i\xi_l\frac{\xi_j\xi_k}{|\xi|^2}e^{-a|\xi|^2}\right],~a>0,~j,k,l=1,2,3.
\end{align*}
Here $\tilde\mu_j(\xi),~j=1,2,3,4,$ are denoted by
\begin{align*}
\mu_j(\xi)&=-\frac{\nu}{4}|\xi|^2+\tilde\mu_j(\xi),~j=1,2, \\
\mu_j(\xi)&=-\frac{\nu+\tilde\nu}{4}|\xi|^2+\tilde\mu_j(\xi),~j=3,4.
\end{align*}
The estimates of $\mathcal{K}_1$ and $\mathcal{K}_2$ follow from Lemma \ref{L1L1LFKS}.
As for $\mathcal{K}_0^{k}$, we use the following $L^1$ estimate of $\mathcal{L}_{l,jk}^a$ shown by Fujigaki and Miyakawa \cite[pp.525--526]{fujigakimiyakawa}.
\begin{lem}\label{fujigakimiyakawa}
Let $a>0$ and $j,k,l=1,2,3$.
Then, the following inequality holds for $t\geq0$:
$$\|\mathcal{L}_{l,jk}^a(\cdot,t)\|_{L^1} \leq C_at^{-\frac{1}{2}}.
$$
\end{lem}
By using Lemma \ref{L1L1LFKS}, Lemma \ref{fujigakimiyakawa} and the Young inequality, we obtain 
\begin{align}
\|\mathcal{F}^{-1}[\hat\varphi_1(\xi)\hat{K}^{23}(\xi,t)\hat\varPsi_0(\xi)]\|_{L^1}&\leq C((1+t)^{\frac{1}{2}}+t^{-\frac{1}{2}})\|\varPsi_0\|_{L^1}, \label{K23L1est1}\\
\|\mathcal{F}^{-1}[\hat\varphi_1(\xi)\hat{K}^{32}(\xi,t)\hat{w}_0(\xi)]\|_{L^1}&\leq C((1+t)^{\frac{1}{2}}+t^{-\frac{1}{2}})\|w_0\|_{L^1}.\label{K32L1est1}
\end{align}
We next show the following uniform bounds with respect to $0\leq t\leq1$:
\begin{align}
\|\mathcal{F}^{-1}[\hat\varphi_1(\xi)\hat{K}^{23}(\xi,t)\hat\varPsi_0(\xi)]\|_{L^1}&\leq C\|\varPsi_0\|_{L^1},~t\geq0, \label{K23bddt}\\
\|\mathcal{F}^{-1}[\hat\varphi_1(\xi)\hat{K}^{32}(\xi,t)\hat{w}_0(\xi)]\|_{L^1}&\leq C\|w_0\|_{L^1},~t\geq0. \label{K32bddt}
\end{align}
To derive \eqref{K23bddt} and \eqref{K32bddt}, we prepare the following lemma proved in \cite{shibatashimizu}.
\begin{lem}\label{shibatashimizu}{\sc (\cite{shibatashimizu})}
Let $\alpha=N+\sigma-3$, where $N\geq0$ is an integer and $0<\sigma\leq1$. Let $f$ be a function such that
\[
f\in C^\infty(\mathbb{R}^3-\{0\}) ,
\]
\[
\partial_\xi^\eta f\in L^1(\mathbb{R}^3),~|\eta|\leq N,
\]
\[
|\partial_\xi^\eta f(\xi)|\leq C_\eta |\xi|^{\alpha-|\eta|},~\xi\neq0.
\]
Then, we have
\[
|\mathcal{F}^{-1}[f(\xi)](x)|\leq C_\alpha \left(\max_{|\eta|\leq N+2} C_\eta\right) |x|^{-3-|\alpha|},~x\neq0.
\]
\end{lem}

By Taylor's formula we have 
\begin{align*}
\frac{e^{\mu_1(\xi)t}-e^{\mu_2(\xi)t}}{\mu_1(\xi)-\mu_2(\xi)}\frac{\xi_j \xi_k\xi_l}{|\xi|^2}\varphi_1(\xi) =\frac{t}{2}e^{-\frac{\nu}{2}|\xi|^2t}\int_0^1e^{i\beta|\xi|f(|\xi|)st}\mathrm{d}s\frac{\xi_j \xi_k\xi_l}{|\xi|^2}\varphi_1(\xi)
\end{align*}
for $|\xi|\leq \frac{M_1}{\sqrt{2}}$, where $f(|\xi|)=\sqrt{1-\frac{\nu^2}{4\beta^2}|\xi|^2}$.

It then follows from the above formula that
\[
\left|\partial_\xi^\eta\left(\frac{e^{\mu_1(\xi)t}-e^{\mu_2(\xi)t}}{\mu_1(\xi)-\mu_2(\xi)}\frac{\xi_j \xi_k\xi_l}{|\xi|^2}\varphi_1(\xi)\right)\right| \leq C_\eta|\xi|^{1-|\eta|}~\mathrm{for}~|\xi|\leq \frac{M_1}{\sqrt{2}}.
\]

We next use Lemma \ref{shibatashimizu} with $(\alpha,N,\sigma)=(1,3,1)$ and calculate in a similar argument as in \cite[pp.228--229]{kobayashishibata} to obtain
\[
\left\|\mathcal{F}^{-1}\left[\frac{e^{\mu_1(\xi)t}-e^{\mu_2(\xi)t}}{\mu_1(\xi)-\mu_2(\xi)}\frac{\xi_j \xi_k\xi_l}{|\xi|^2}\varphi_1(\xi)  \right]\right\|_{L^1}\leq C,~0\leq t\leq 1.
\]
Similarly, we can prove
\[
\left\|\mathcal{F}^{-1}\left[\frac{e^{\mu_3(\xi)t}-e^{\mu_4(\xi)t}}{\mu_3(\xi)-\mu_4(\xi)}\frac{\xi_j \xi_k\xi_l}{|\xi|^2}\varphi_1(\xi)  \right]\right\|_{L^1}\leq C,~0\leq t\leq 1.
\]
We thus arrive at \eqref{K23bddt} and \eqref{K32bddt}.

By \eqref{K23L1est1}--\eqref{K32bddt}, we have
\begin{align}
\|\mathcal{F}^{-1}[\hat\varphi_1(\xi)\hat{K}^{23}(\xi,t)\hat\varPsi_0(\xi)]\|_{L^1}&\leq C(1+t)^{\frac{1}{2}}\|\varPsi_0\|_{L^1},~t\geq0, \label{K32L1est2}\\
\|\mathcal{F}^{-1}[\hat\varphi_1(\xi)\hat{K}^{32}(\xi,t)\hat{w}_0(\xi)]\|_{L^1}&\leq C(1+t)^{\frac{1}{2}}\|w_0\|_{L^1},~t\geq0. \label{K23L1est2}
\end{align}

We see from \eqref{ineq6.7}--\eqref{ineq6.11}, \eqref{K32L1est2} and \eqref{K23L1est2} that
\begin{align*}
\|e^{-tL}U_1(0)\|_{L^1}\leq C(1+t)^{\frac{1}{2}}\|u_0\|_{L^1}.
\end{align*}
This completes the proof of Lemma \ref{L1L1est-low-semigroup}. $\blacksquare$

We next estimate $\int_{0}^{t} \|e^{-(t-s)L}P_1N(s)\|_{L^1}\mathrm{d}s$.

\begin{lem}\label{estlowfreqP1NL1}
There exists a positive number $\delta_0$ such that if $\|u_0\|_{\mathcal{X}}\leq\delta_0$, then the following estimate holds:
\[
\int_{0}^{t} \|e^{-(t-s)L}P_1N(s)\|_{L^1}\mathrm{d}s\leq C(1+t)^{\frac{1}{2}}\|u_0\|_{\mathcal{X}},~t\geq0.
\]
\end{lem}

{\bf Proof.} We obtain the following estimate in a similar argument as in the proof of Lemma \ref{NLestLFLinf} by using \eqref{calK23}, \eqref{calK33}, \eqref{ineq5.12}, \eqref{ineq5.13} and  Lemma \ref{L1L1LFKS} : 
\begin{align}
\begin{array}{l}
\displaystyle
\left\|\mathcal{F}^{-1}\left[\hat\varphi_1(\xi)\hat{K}^{jk}(\xi,t-s)\hat{N}_k(\xi,s)\right]\right\|_{L^1}  \\[1ex]
\quad \leq C(1+t-s)^{-2}(1+s)^{\frac{1}{2}}\|u_0\|_{\mathcal{X}},~j,k=1,2,3.
\end{array}
\label{estN111L} 
\end{align}
By using Lemma \ref{betafunc1} with $a=b=2$, we have 
\begin{align}
\begin{array}{l}
\displaystyle
\int_0^t (1+t-s)^{-2}(1+s)^{\frac{1}{2}}\mathrm{d}s \\[2ex]
\quad\quad
\displaystyle
\leq
(1+t)^{\frac{5}{2}}\int_0^t (1+t-s)^{-2}(1+s)^{-2}\mathrm{d}s \\[2ex]
\quad\quad
\displaystyle
\leq C(1+t)^{\frac{1}{2}}. 
\end{array}
\label{ineqLEM671}
\end{align}
We then see from \eqref{estN111L} and \eqref{ineqLEM671} that  
\[
\int_{0}^{t} \|e^{-(t-s)L}P_1N(s)\|_{L^1}\mathrm{d}s\leq C(1+t)^{\frac{1}{2}}\|u_0\|_{\mathcal{X}},~t\geq0.
\]
This completes the proof. $\blacksquare$

\textbf{Proof of Proposition \ref{L1Lplowfreq}.} Taking $L^1$ norm of the first equation of \eqref{lowfreq}, we obtain
\begin{align*}
\|U_1(t)\|_{L^1}&\leq \|e^{-tL}U_1(0)\|_{L^1}+\int_{0}^{t} \|e^{-(t-s)L}P_1N(s)\|_{L^1}\mathrm{d}s. 
\end{align*}
This completes the proof. $\blacksquare$

We next consider the high frequency part $U_\infty$.
\begin{prop}\label{LpLphighfreq}
There exists a positive number $\delta_p$ such that if $\|u_0\|_{\mathcal{X}}\leq\delta_p$, then it holds the following estimate for $t\geq0$:
\begin{align*}
\|U_\infty(t)\|_{L^p}\leq C(1+t)^{-\frac{3}{2}\left(1-\frac{1}{p}\right)+\frac{1}{2}\left(\frac{2}{p}-1\right)}(\|u_0\|_{L^p}+\|u_0\|_{\mathcal{X}}).
\end{align*}
\end{prop}
In order to prove Proposition \ref{LpLphighfreq}, we prepare the following lemma.
\begin{lem}\label{highfreqshibata2000}
Let $1<p<\infty$ and $f\in L^p$. Then, the following estimates hold for $t\geq0$:
\begin{align*}
&\left\|\partial_t^j\partial_x^\alpha\mathcal{F}^{-1}\left[\frac{e^{\mu_1(\xi)t}-e^{\mu_2(\xi)t}}{\mu_1(\xi)-\mu_2(\xi)} \hat\varphi_\infty(\xi)\hat{f}(\xi)\right]\right\|_{L^p} \leq Ce^{-ct}\|f\|_{L^p}~j+|\alpha|=1, 
\\ 
&\left\|\mathcal{F}^{-1}\left[\frac{\mu_1(\xi)e^{\mu_2(\xi)t}-\mu_2(\xi)e^{\mu_1(\xi)t}}{\mu_1(\xi)-\mu_2(\xi)}\hat\varphi_\infty(\xi)\hat{f}(\xi)\right]\right\|_{L^p} \leq Ce^{-ct}\|f\|_{L^p}, 
\\ 
&\left\|\partial_t^j\partial_x^\alpha\mathcal{F}^{-1}\left[\frac{e^{\mu_3(\xi)t}-e^{\mu_4(\xi)t}}{\mu_3(\xi)-\mu_4(\xi)} \hat\varphi_\infty(\xi)\hat{f}(\xi)\right]\right\|_{L^p} \leq Ce^{-ct}\|f\|_{L^p},~j+|\alpha|=1, 
\\ 
&\left\|\mathcal{F}^{-1}\left[\frac{\mu_3(\xi)e^{\mu_4(\xi)t}-\mu_4(\xi)e^{\mu_3(\xi)t}}{\mu_3(\xi)-\mu_4(\xi)}\hat\varphi_\infty(\xi)\hat{f}(\xi)\right]\right\|_{L^p} \leq Ce^{-ct}\|f\|_{L^p}. 
\end{align*}
\end{lem}
Lemma \ref{highfreqshibata2000} directly follows from \cite[Theorem 4.1]{shibataVW2000}.
We first consider $\|e^{-tL}U_\infty(0)\|_{L^p}$.
\begin{lem}\label{LpLpest-high-semigroup}
The following estimate holds for $t\geq0$:
\begin{align}
\|e^{-tL} U_\infty(0)\|_{L^p}\leq Ce^{-ct}\|u_0\|_{L^p}. \label{SGHFLpLp}
\end{align}
\end{lem}
\textbf{Proof.} The estimate \eqref{SGHFLpLp} can be shown by using Lemma \ref{psitildepsi}, Lemma \ref{highfreqshibata2000} and the $L^p$ boundedness of the Riesz operator.
This completes the proof. $\blacksquare$

We next estimate $\int_{0}^{t} \|e^{-(t-s)L}P_\infty N(s)\|_{L^p}\mathrm{d}s$.

\begin{lem}\label{estlowfreqPinfNLp}
There exists a positive number $\delta_p$ such that if $\|u_0\|_{\mathcal{X}}\leq\delta_p$, then the following estimate holds:
\[
\int_{0}^{t} \|e^{-(t-s)L}P_\infty N(s)\|_{L^p}\mathrm{d}s\leq C(1+t)^{-\frac{3}{2}\left(1-\frac{1}{p}\right)+\frac{1}{2}\left(\frac{2}{p}-1\right)}\|u_0\|_{\mathcal{X}},~t\geq0.
\]
\end{lem}
{\bf Proof.} We obtain the following estimate in a similar argument as in the proof of Lemma \ref{LpLpest-high-semigroup}: 
\begin{align}
\begin{array}{l}
\displaystyle
\left\|\mathcal{F}^{-1}\left[\hat\varphi_\infty(\xi)\hat{K}^{jk}(\xi,t-s)\hat{N}_k(\xi,s)\right]\right\|_{L^p} \leq Ce^{-c(t-s)}\|N_k(s)\|_{L^p}, \\[2ex]
\hspace{7cm} j=1,2,3,~k=1,2. 
\end{array}
\label{ineqLEM6110} 
\end{align}
In view of Lemma \ref{Sobolevineq} and the $L^p$ boundedness of the Riesz operator, we have
\begin{align}
\|N_k(s)\|_{L^p}\leq C\|U(s)\|_{H^2}^2\leq C(1+t)^{-\frac{3}{2}}\|u_0\|_{\mathcal{X}},~k=1,2,3. \label{ineqLEM6111}
\end{align}
By employing Lemma \ref{betafunc1} with $a=b=\frac{3}{2}$, we have 
\begin{align}
\begin{array}{l}
\displaystyle
\int_0^t e^{-c(t-s)}(1+s)^{-\frac{3}{2}}\mathrm{d}s
\leq C(1+t)^{-\frac{3}{2}} \\[2ex]
\quad\quad\quad\quad\quad\quad\quad\quad\quad\quad
\displaystyle
\leq C(1+t)^{-2+\frac{5}{2p}} \\[2ex]
\quad\quad\quad\quad\quad\quad\quad\quad\quad\quad
\displaystyle
=C(1+t)^{-\frac{3}{2}\left(1-\frac{1}{p}\right)+\frac{1}{2}\left(\frac{2}{p}-1\right)}.
\end{array}
\label{ineqLEM6112}
\end{align}
Together with \eqref{ineqLEM6110}--\eqref{ineqLEM6112} yields
\begin{align*}
\int_{0}^{t} \|e^{-(t-s)L}P_\infty N(s)\|_{L^p}\mathrm{d}s\leq C(1+s)^{-\frac{3}{2}\left(1-\frac{1}{p}\right)+\frac{1}{2}\left(\frac{2}{p}-1\right)}\|u_0\|_{\mathcal{X}}.
\end{align*}
This completes the proof. $\blacksquare$

\textbf{Proof of Proposition \ref{LpLphighfreq}.}
By taking $L^p$ norm of the first equation of \eqref{highfreq}, we have
\begin{align}\label{lowfreqineqL1}
\|U_\infty(t)\|_{L^p}&\leq \|e^{-tL}U_\infty(0)\|_{L^p}+\int_{0}^{t} \|e^{-(t-s)L}P_\infty N(s)\|_{L^p}\mathrm{d}s. 
\end{align}
Combining Lemma \ref{LpLpest-high-semigroup}, Lemma \ref{estlowfreqPinfNLp} and \eqref{lowfreqineqL1}, we arrive at
\[
\|U_\infty(t)\|_{L^p}\leq C(1+t)^{-\frac{3}{2}\left(1-\frac{1}{p}\right)+\frac{1}{2}\left(\frac{2}{p}-1\right)}(\|u_0\|_{L^p}+\|u_0\|_{\mathcal{X}}),~t\geq0.
\]
This completes the proof of Proposition \ref{LpLphighfreq}. $\blacksquare$
\\[2ex]
{\large \bf Appendix: Proof of Lemma \ref{solform}} \\[2ex]
In this section, we derive the solution formula \eqref{solformula}. \\
\textbf{Proof of Lemma \ref{solform}.}
We write \eqref{Flinearizedproblem5} as
\begin{gather}
\partial_t\hat{\phi}+i\xi\cdot \hat{w}=0,\tag{A.1}
\label{AP1}\\[1ex]
\partial_t  \hat{w}+\nu|\xi|^2 \hat{w}
+\tilde \nu\xi{}^\top\xi  \hat{w}
+i\gamma^2\xi\hat{\phi}
-i\beta^2 \hat{\varPsi}\xi  = 0,\tag{A.2} 
\label{AP2}\\[1ex]
\partial_t \hat{\varPsi} -i\hat{w}^\top\xi= 0, \tag{A.3}
\label{AP3}\\[1ex]
\hat{\phi}+i\xi\cdot \hat{\psi}=0,~\hat{\varPsi}=i\hat{\psi}^\top\xi, \tag{A.4}
\label{AP4}\\[1ex]
(\hat\phi, \hat{w}, \hat\varPsi)|_{t=0}=(\hat\phi_0, \hat{w}_0, \hat\varPsi_0),~\hat{\phi}_0+i\xi\cdot \hat{\psi}_0=0.\tag{A.5}
\label{AP5}
\end{gather}
Setting $w_t=\partial_t w$, we see from \eqref{AP2}--\eqref{AP4} that
\begin{gather}
\partial_t
\left(
\begin{array}{l}
\hat{w}\\
\hat{w}_t
\end{array}
\right)+\mathcal{A}(\xi)
\left(
\begin{array}{l}
\hat{w}\\
\hat{w}_t
\end{array}
\right)
=0,
~\left.\left(
\begin{array}{l}
\hat{w}\\
\hat{w}_t
\end{array}
\right)\right|_{t=0}
=\left(
\begin{array}{l}
\hat{w_0}\\
\hat{w}_{t,0}
\end{array}
\right). \tag{A.6} \label{AP2.1}
\end{gather}
Here
\[
\mathcal{A}(\xi)=\left(
\begin{array}{@{\ }cc@{\ } }
0 & -I\\
\beta^2|\xi|^2I+\gamma^2\xi{}^\top\xi & \nu|\xi|^2I+\tilde\nu\xi{}^\top\xi
\end{array}
\right)
\]
and
\begin{gather}
\hat{w}_{t,0}=
-i\gamma^2\xi\hat{\phi}_0-(\nu|\xi|^2I 
+\tilde \nu\xi{}^\top\xi  )\hat{w}_0
+i\beta^2 \hat{\varPsi}_0\xi.
\tag{A.7}\label{AP8}
\end{gather}
To solve \eqref{solformula}, we first investigate the characteristic equation of $-\mathcal{A}(\xi)$. Let $T$ be a $3\times3$ orthogonal matrix and set
\[
\mathcal{T}=\left(
\begin{array}{@{\ }cc@{\ } }
T & 0\\
0 & T
\end{array}
\right).
\]
We see that
\[
\mathcal{A}(T\xi)=\mathcal{T}\mathcal{A}(\xi){}^\top\mathcal{T}.
\] 
We choose $T$ so that $T\xi=r e_1$, where $r=|\xi|$ and $e_1={}^\top(1,0,0)$. Using this $T$, we have
\begin{align*}
\det(\mu I_6+\mathcal{A}(\xi)) &=\det(\mathcal{T}(\mu I_6+\mathcal{A}(\xi)){}^\top\mathcal{T})\\
&=\det(\mu I_6+\mathcal{A}(T\xi)) \\
&=\det\left(
\begin{array}{@{\ }cc@{\ } }
\mu I & -I\\
\beta^2r^2I+\gamma^2r^2e_1{}^\top e_1 & (\mu +\nu r^2)I+\tilde\nu r^2e_1{}^\top e_1
\end{array}
\right) \\
&=(\mu^2+\nu r^2\mu+\beta^2r^2)^2[\mu^2+(\nu+\tilde\nu) r^2\mu+(\beta^2+\gamma^2)r^2]\\
&=(\mu-\mu_1(\xi))^2(\mu-\mu_2(\xi))^2(\mu-\mu_3(\xi))(\mu-\mu_4(\xi)).
\end{align*}  
Therefore, the eigenvalues of $-\mathcal{A}(\xi)$ are given by $\mu_j(\xi)$, $j=1,2,3,4$. We note that
\begin{align}
&\mu_1(\xi)\mu_2(\xi)=\beta^2|\xi|^2,~
\mu_1(\xi)+\mu_2(\xi)=-\nu|\xi|^2, \tag{A.8} \label{AP2.2}\\
&\mu_3(\xi)\mu_4(\xi)=(\beta^2+\gamma^2)|\xi|^2,~ 
\mu_3(\xi)+\mu_4(\xi)=-(\nu+\tilde{\nu})|\xi|^2. \tag{A.9} \label{AP2.3}
\end{align}
By using \eqref{AP2.2} and \eqref{AP2.3}, the eigenprojections for $\mu_j(\xi)$ of $-\mathcal{A}(\xi)$ are written by
\[
\Pi_1(\xi)=\frac{1}{\mu_1(\xi)-\mu_2(\xi)}
\left(
\begin{array}{@{\ }cc@{\ } }
\displaystyle
-\mu_2(\xi)\left( I-\frac{\xi{}^\top\xi}{|\xi|^2}\right) &\displaystyle I-\frac{\xi{}^\top\xi}{|\xi|^2}\\[1ex]
\displaystyle
-\mu_1(\xi)\mu_2(\xi)\left( I-\frac{\xi{}^\top\xi}{|\xi|^2}\right) &\displaystyle \mu_1(\xi)\left(I-\frac{\xi{}^\top\xi}{|\xi|^2}\right)
\end{array}
\right),
\]
\[
\Pi_2(\xi)=\frac{1}{\mu_1(\xi)-\mu_2(\xi)}
\left(
\begin{array}{@{\ }cc@{\ } }
\displaystyle
\mu_1(\xi)\left( I-\frac{\xi{}^\top\xi}{|\xi|^2}\right) &\displaystyle-\left( I-\frac{\xi{}^\top\xi}{|\xi|^2}\right)\\[1ex]
\displaystyle
\mu_1(\xi)\mu_2(\xi)\left( I-\frac{\xi{}^\top\xi}{|\xi|^2}\right) &\displaystyle -\mu_2(\xi)\left(I-\frac{\xi{}^\top\xi}{|\xi|^2}\right)
\end{array}
\right),
\]
\[
\Pi_3(\xi)=\frac{1}{\mu_3(\xi)-\mu_4(\xi)}
\left(
\begin{array}{@{\ }cc@{\ } }
\displaystyle
-\mu_4(\xi)\frac{\xi{}^\top\xi}{|\xi|^2} &\displaystyle \frac{\xi{}^\top\xi}{|\xi|^2}\\[1ex]
\displaystyle
-\mu_3(\xi)\mu_4(\xi)\frac{\xi{}^\top\xi}{|\xi|^2} &\displaystyle \mu_3(\xi)\frac{\xi{}^\top\xi}{|\xi|^2}
\end{array}
\right),
\]
\[
\Pi_4(\xi)=\frac{1}{\mu_3(\xi)-\mu_4(\xi)}
\left(
\begin{array}{@{\ }cc@{\ } }
\displaystyle
\mu_3(\xi)\frac{\xi{}^\top\xi}{|\xi|^2} & \displaystyle -\frac{\xi{}^\top\xi}{|\xi|^2}\\[1ex]
\displaystyle
\mu_3(\xi)\mu_4(\xi)\frac{\xi{}^\top\xi}{|\xi|^2} &\displaystyle -\mu_4(\xi)\frac{\xi{}^\top\xi}{|\xi|^2}
\end{array}
\right).
\]
The solution semigroup $e^{-t\Hat{\mathcal{A}}(\xi) }$ is then expressed as
\begin{align*}
&e^{-t\Hat{\mathcal{A}}(\xi) }\\[1ex]
&=e^{\mu_1(\xi)t}\Pi_1(\xi)+e^{\mu_2(\xi)t}\Pi_2(\xi)+e^{\mu_3(\xi)t}\Pi_3(\xi)+e^{\mu_4(\xi)t}\Pi_4(\xi) \\[1ex]
                     &=\left(
\begin{array}{@{\ }cc@{\ } }
\displaystyle
\frac{\mu_1(\xi)e^{\mu_2(\xi)t}-\mu_2(\xi)e^{\mu_1(\xi)t}}{\mu_1(\xi)-\mu_2(\xi)}\left( I-\frac{\xi{}^\top\xi}{|\xi|^2}\right) &\displaystyle\frac{e^{\mu_1(\xi)t}-e^{\mu_2(\xi)t}}{\mu_1(\xi)-\mu_2(\xi)}\left( I-\frac{\xi{}^\top\xi}{|\xi|^2}\right)\\[2ex]
\displaystyle
-\mu_1(\xi)\mu_2(\xi)\frac{e^{\mu_1(\xi)t}-e^{\mu_2(\xi)t}}{\mu_1(\xi)-\mu_2(\xi)}\left( I-\frac{\xi{}^\top\xi}{|\xi|^2}\right) &\displaystyle \frac{\mu_1(\xi)e^{\mu_1(\xi)t}-\mu_2(\xi)e^{\mu_2(\xi)t}}{\mu_1(\xi)-\mu_2(\xi)}\left(I-\frac{\xi{}^\top\xi}{|\xi|^2}\right)
\end{array}
\right) \\[1ex]
&\quad+\left(
\begin{array}{@{\ }cc@{\ } }
\displaystyle
\frac{\mu_3(\xi)e^{\mu_4(\xi)t}-\mu_4(\xi)e^{\mu_3(\xi)t}}{\mu_3(\xi)-\mu_4(\xi)}\frac{\xi{}^\top\xi}{|\xi|^2} &\displaystyle \frac{e^{\mu_3(\xi)t}-e^{\mu_4(\xi)t}}{\mu_3(\xi)-\mu_4(\xi)}\frac{\xi{}^\top\xi}{|\xi|^2}\\[2ex]
\displaystyle
-\mu_3(\xi)\mu_4(\xi)\frac{e^{\mu_3(\xi)t}-e^{\mu_4(\xi)t}}{\mu_3(\xi)-\mu_4(\xi)}\frac{\xi{}^\top\xi}{|\xi|^2} &\displaystyle \frac{\mu_3(\xi)e^{\mu_3(\xi)t}-\mu_4(\xi)e^{\mu_4(\xi)t}}{\mu_3(\xi)-\mu_4(\xi)}\frac{\xi{}^\top\xi}{|\xi|^2}
\end{array}
\right)
\end{align*}
It then follows that $\hat{w}(\xi,t)$ is written as
\begin{gather}
\begin{array}{l}
\displaystyle
\hat{w}(\xi,t)=\frac{\mu_1(\xi)e^{\mu_2(\xi)t}-\mu_2(\xi)e^{\mu_1(\xi)t}}{\mu_1(\xi)-\mu_2(\xi)}\left( I-\frac{\xi{}^\top\xi}{|\xi|^2}\right)\hat{w}_0(\xi)\\[2ex]
\displaystyle
\quad\quad\quad\quad
+\frac{e^{\mu_1(\xi)t}-e^{\mu_2(\xi)t}}{\mu_1(\xi)-\mu_2(\xi)}\left( I-\frac{\xi{}^\top\xi}{|\xi|^2}\right)\hat{w}_{t,0}(\xi)\\[2ex]
\displaystyle
\quad\quad\quad\quad\quad
+\frac{\mu_3(\xi)e^{\mu_4(\xi)t}-\mu_4(\xi)e^{\mu_3(\xi)t}}{\mu_3(\xi)-\mu_4(\xi)}\frac{\xi{}^\top\xi}{|\xi|^2}\hat{w}_0(\xi)\\[2ex]
\displaystyle
\quad\quad\quad\quad\quad\quad+\frac{e^{\mu_3(\xi)t}-e^{\mu_4(\xi)t}}{\mu_3(\xi)-\mu_4(\xi)}\frac{\xi{}^\top\xi}{|\xi|^2}\hat{w}_{t,0}(\xi).
\end{array}
\tag{A.10}\label{AP7}
\end{gather}
Substituting \eqref{AP8} into \eqref{AP7} leads to
\begin{gather}
\begin{array}{l}
\displaystyle
\hat{w}(\xi,t)=-i\gamma^2\frac{e^{\mu_3(\xi)t}-e^{\mu_4(\xi)t}}{\mu_3(\xi)-\mu_4(\xi)}\xi\hat{\phi}_0(\xi)\\[2ex]
\displaystyle
\quad\quad\quad+\frac{\mu_1(\xi)e^{\mu_1(\xi)t}-\mu_2(\xi)e^{\mu_2(\xi)t}}{\mu_1(\xi)-\mu_2(\xi)}\left(I-\frac{{\xi}^\top \xi}{|\xi|^2}\right)\hat{w}_0(\xi) \\[2ex]
\displaystyle
\quad\quad\quad\quad+\frac{\mu_3(\xi)e^{\mu_3(\xi)t}-\mu_4(\xi)e^{\mu_4(\xi)t}}{\mu_3(\xi)-\mu_4(\xi)}\frac{{\xi}^\top \xi}{|\xi|^2}\hat{w}_0(\xi)\\[2ex]
\displaystyle
\quad\quad\quad+i\beta^2\frac{e^{\mu_1(\xi)t}-e^{\mu_2(\xi)t}}{\mu_1(\xi)-\mu_2(\xi)}\left(I-\frac{{\xi}^\top \xi}{|\xi|^2}\right)\hat{\varPsi}_0(\xi)\xi \\[2ex]
\displaystyle
\quad\quad\quad\quad+i\beta^2\frac{e^{\mu_3(\xi)t}-e^{\mu_4(\xi)t}}{\mu_3(\xi)-\mu_4(\xi)}\frac{{\xi}^\top \xi}{|\xi|^2}\hat{\varPsi}_0(\xi)\xi.
\end{array}
\tag{A.11}\label{AP9}
\end{gather}
We see from \eqref{AP1}, \eqref{AP3}, \eqref{AP5} and \eqref{AP9} that 
\begin{align*}
\hat\phi(\xi,t)&=\hat\phi_0(\xi)-i\xi\cdot\int_0^t \hat{w}(\xi,s) \mathrm{d}s  \\
&=\frac{\beta^2}{\beta^2+\gamma^2}\left(\hat\phi_0(\xi)+\frac{{}^\top\xi\hat\varPsi_0(\xi)\xi}{|\xi|^2}\right) \\
&\quad+\frac{\mu_3(\xi)e^{\mu_4(\xi)t}-\mu_4(\xi)e^{\mu_3(\xi)t}}{\mu_3(\xi)-\mu_4(\xi)}\left(\frac{\gamma^2}{\beta^2+\gamma^2}\hat\phi_0(\xi)-\frac{\beta^2}{\beta^2+\gamma^2}\frac{{}^\top\xi\hat\varPsi_0(\xi)\xi}{|\xi|^2}\right)  \\
&\quad\quad-i\frac{e^{\mu_3(\xi)t}-e^{\mu_4(\xi)t}}{\mu_3(\xi)-\mu_4(\xi)}\xi\cdot\hat{w}_0(\xi)\\
&=\frac{\mu_3(\xi)e^{\mu_4(\xi)t}-\mu_4(\xi)e^{\mu_3(\xi)t}}{\mu_3(\xi)-\mu_4(\xi)}\hat\phi_0(\xi) 
-i\frac{e^{\mu_3(\xi)t}-e^{\mu_4(\xi)t}}{\mu_3(\xi)-\mu_4(\xi)}\xi\cdot\hat{w}_0(\xi),
\end{align*}
\begin{align*}
\hat\varPsi(\xi,t)&=\hat\varPsi_0(\xi)+i\left(\int_0^t \hat{w}(\xi,s) \mathrm{d}s\right){}^\top\xi  \\
&=\frac{\gamma^2}{\beta^2+\gamma^2}\left(\hat\phi_0(\xi)\frac{\xi{}^\top\xi}{|\xi|^2}+\frac{\xi{}^\top\xi}{|\xi|^2}\hat\varPsi_0(\xi)\right) \\
&\quad+i\frac{e^{\mu_1(\xi)t}-e^{\mu_2(\xi)t}}{\mu_1(\xi)-\mu_2(\xi)}\left(I-\frac{{\xi}^\top \xi}{|\xi|^2}\right)\hat{w}_0(\xi){}^\top\xi \\
&\quad\quad+i\frac{e^{\mu_3(\xi)t}-e^{\mu_4(\xi)t}}{\mu_3(\xi)-\mu_4(\xi)}\frac{{\xi}^\top \xi}{|\xi|^2}\hat{w}_0(\xi){}^\top\xi\\
&\quad+\frac{\mu_1(\xi)e^{\mu_2(\xi)t}-\mu_2(\xi)e^{\mu_1(\xi)t}}{\mu_1(\xi)-\mu_2(\xi)}\left(I-\frac{{\xi}^\top \xi}{|\xi|^2}\right)\hat\varPsi_0(\xi)\\
&\quad\quad+\frac{\mu_3(\xi)e^{\mu_4(\xi)t}-\mu_4(\xi)e^{\mu_3(\xi)t}}{\mu_3(\xi)-\mu_4(\xi)}\left(-\frac{\gamma^2}{\beta^2+\gamma^2}\hat\phi_0(\xi)\frac{\xi{}^\top\xi}{|\xi|^2}+\frac{\beta^2}{\beta^2+\gamma^2}\frac{{\xi}^\top \xi}{|\xi|^2}\hat\varPsi_0(\xi)\right)\\
&=i\frac{e^{\mu_1(\xi)t}-e^{\mu_2(\xi)t}}{\mu_1(\xi)-\mu_2(\xi)}\left(I-\frac{{\xi}^\top \xi}{|\xi|^2}\right)\hat{w}_0(\xi){}^\top\xi \\
&\quad\quad+i\frac{e^{\mu_3(\xi)t}-e^{\mu_4(\xi)t}}{\mu_3(\xi)-\mu_4(\xi)}\frac{{\xi}^\top \xi}{|\xi|^2}\hat{w}_0(\xi){}^\top\xi\\\quad
&\quad+\frac{\mu_1(\xi)e^{\mu_2(\xi)t}-\mu_2(\xi)e^{\mu_1(\xi)t}}{\mu_1(\xi)-\mu_2(\xi)}\left(I-\frac{{\xi}^\top \xi}{|\xi|^2}\right)\hat\varPsi_0(\xi)\\
&\quad\quad+\frac{\mu_3(\xi)e^{\mu_4(\xi)t}-\mu_4(\xi)e^{\mu_3(\xi)t}}{\mu_3(\xi)-\mu_4(\xi)}\frac{{\xi}^\top \xi}{|\xi|^2}\hat\varPsi_0(\xi).
\end{align*}
This completes the proof. $\blacksquare$ \\[2ex]
\textbf{Acknowledgements. } The author would like to thank Professor Yoshiyuki Kagei for his valuable suggestions and comments. 
This work was partially supported by JSPS KAKENHI Grant Number 19J10056. 


\begin{thebibliography}{9}
\bibitem{fujigakimiyakawa} Y. Fujigaki, T. Miyakawa, Asymptotic profiles of nonstationary incompressible Navier-Stokes flows in the whole space, SIAM J. Math. Anal., $\bf{33}$ (2001), pp. 523--544.
\bibitem{gurtin}M. E. Gurtin, \textit{An Introduction to Continuum Mechanics}, Math. Sci. Eng., vol. 158, Academic Press, New York--London, 1981.
\bibitem{hoffzumbrun}  D. Hoff, K. Zumbrun, Multi-dimensional diffusion waves for the Navier-Stokes equations of compressible ﬂow, Indiana Univ. Math. J., $\bf{44}$ (1995), pp. 603--676.
\bibitem{hu} X. Hu, Global existence of weak solutions to two dimensional compressible viscoelastic flows, J. Differential Equations, $\bf{265}$ (2018), pp. 3130--3167.
\bibitem{huwang1} X. Hu, D. Wang, Local strong solution to the compressible viscoelastic flow with large data, J. Differential Equations, $\bf{249}$ (2010), pp. 1179--1198.
\bibitem{huwang2} X. Hu, D. Wang, Global existence for the multi-dimensional compressible viscoelastic flows, J. Differential Equations, $\bf{250}$ (2011), pp. 1200--1231.
\bibitem{huwu}X. Hu, G. Wu, Global existence and optimal decay rates for three-dimensional compressible viscoelastic flows, SIAM J. Math. Anal., $\bf{45}$~(2013),  pp. 2815--2833.
\bibitem{kawashimamatsumuranishida} S. Kawashima, A. Matsumura, T. Nishida, On the fluid-dynamical approximation to the Boltzmann equation at the level of the Navier-Stokes equation, Commun. Math. Phys., $\bf{70}$, (1979), pp.97--124.  
\bibitem{kobayashishibata}T. Kobayashi, Y. Shibata, Remark on the rate of decay of solutions to linearized compressible Navier-Stokes equation, Pacific J. Math. $\bf{207}$ (2002), pp. 199--234.
\bibitem{liweiyao}Y. Li, R. Wei, Z. Yao, Optimal decay rates for the compressible viscoelastic flows, J. Math. Phys., $\bf{57}$, 111506, (2016).
\bibitem{linliuzhang}F. Lin, C. Liu, P. Zhang, On hydrodynamics of viscoelastic fluids, Comm. Pure Appl. Math., $\bf{58}$ (2005), pp. 1437--1471.
\bibitem{liuwuxuzhang}L. Liu, Y. Wu, F. Xu, X. Zhang, The optimal convergence rates for the multi-dimensional compressible viscoelastic ﬂows, Z. Angew. Math. Mech. $\bf{96}$, (2016), pp.1490–-1504.
\bibitem{matsumuranishida1979}A. Matsumura, T. Nishida, The initial value problems for the equation of motion of compressible viscous and heat-conductive fluids, Proc. Japan Acad. Ser., {\bf 89}~(1979), pp. 337--342.
\bibitem{matsumura-nishidaEM}A. Matsumura, T. Nishida, Initial value problem for the equations of motion of  viscous and heat-conductive gases, J. Math. Kyoto Univ., {\bf 20}~(1980), pp. 67--104.
\bibitem{panxu} X. Pan, J. Xu, Global existence and optimal decay estimates of the compressible viscoelastic flows in $L^p$ critical spaces, Discrete Contin. Dyn. Syst., $\bf{39}$, (2019), pp. 2021--2057.
\bibitem{qianzhang}J. Qian, Z. Zhang, Global well-posedness for compressible viscoelastic fluids near equilibrium, Arch. Ration. Mech. Anal., $\bf{198}$ (2010), pp. 835--868. 
\bibitem{qian}J. Qian, Initial boundary value problems for the compressible viscoelastic fluid, J. Differential Equations, $\bf{250}$ (2011), pp. 848--865. 
\bibitem{shibataVW2000} Y. Shibata, On the rate of decay of solutions to linear viscoelastic equation, Math. Methods Appl. Sci., $\bf{23}$ (2000), pp. 203--226. 
\bibitem{shibatashimizu} Y. Shibata, Y. Shimizu, A decay property of the Fourier transform and its application to Stokes problem, J. Math. Fluid Mech., $\bf{3}$ (2001), pp. 213--230.
\bibitem{sideristhomases}T. C. Sideris, B. Thomases, Global existence for 3D incompressible isotropic elastodynamics via the incompressible limit, Comm. Pure Appl. Math., $\bf{57}$~(2004), pp. 1--39.
\bibitem{segal}  I. Segal, Dispersion for non-linear relativistic equations. II, Ann. Sci. ´Ecole Norm. Sup. (4), $\bf{1}$ (1968), pp. 459–-497. 
\bibitem{weiliyao} R. Wei, Y. Li, Z. Yao, Decay of the compressible viscoelastic flows, Commun. Pure Appl. Anal., $\bf{15}$ (2016), pp. 1603–-1624.
\bibitem{wugaotan} G. Wu, Z. Gao, Z. Tan, Time decay rates for the compressible viscoelastic flows, J. Math. Anal. Appl., $\bf{452}$ (2017),  pp. 990–-1004.
\end{thebibliography}
\end{document}